\documentclass[12pt]{amsart}

\usepackage{latexsym}
\usepackage{amsmath}
\usepackage{amssymb}
\usepackage{amscd}
\usepackage{array}
\usepackage{hhline}
\usepackage{color}
\usepackage[all]{xy}
\usepackage{enumerate}
\usepackage{graphicx}
\usepackage{lscape}

\begin{document}

\title{the triplet vertex operator algebra $W(p)$ and the restricted quantum group $\bar{U}_q(sl_2)$ at $q=e^{\frac{\pi i}{p}}$}
\author{Kiyokazu Nagatomo and Akihiro Tsuchiya} 
\date{}
\address{Graduate School of Information Science and Technology, Osaka University, Toyonaka, Osaka 560-0043, Japan}
\email{nagatomo@math.sci.osaka-u.ac.jp}
\address{Institute for the Physics and Mathematics of the Universe, University of Tokyo, Kashiwa City, Chiba 277-8568, Japan}
\email{akihiro.tsuchiya@ipmu.jp}

\maketitle
\begin{abstract}
We study the abelian category $W(p)$-mod of modules over the triplet $W$ algebra $W(p)$.  We construct the projective covers $\mathcal{P}_s^\pm$ of all the simple objects $\mathcal{X}_s^\pm$, $1\leqq s\leqq p$, in the category $W(p)$-mod.  By using the structure of these projective modules, we show that $W(p)$-mod is a category which is equivalent to the abelian category of the finite-dimensional modules for the restricted quantum group $\bar U_q(sl_2)$ at $q=e^{\frac{\pi i}{p}}$.  This Kazdan-Lusztig type correspondence was conjectured by Feigin et al. [FGST1], [FGST2].
\end{abstract}
\vspace{0.2in}
\section{Introduction}
The theory of vertex operator algebra (VOA) is an algebraic counter-part of comformal field theory.  About general facts around VOA, see [FrB].  Up to now, examples of conformal field theory over general Riemann surfaces are constructed by using lattice VOAs, VOAs associated with integrable representations of affine Lie algebras with the positive integer level, or VOAs associated with the minimal series of the Virasoro algebra.  The abelian category of modules over these VOA's are all semi-simple and the number of simple objects is finite.  In order to define a conformal field theory on Riemann surfaces associated with a VOA, it is necessary that this VOA has some finiteness condition.  Zhu found such a finiteness condition on a VOA called the $C_2$-finiteness condition, and showed that the abelian category of modules over a VOA satisfying $C_2$-finiteness condition is Artinian and Noetherian, moreover, the number of simple objects is finite [FrZ], [Zhu].
\par
Associated to a VOA which has $C_2$-finiteness condition, Zhu developed the theory of conformal blocks on Riemann surfaces, and showed that the dimension of conformal blocks are finite for genus one case, and found the K-Z type differential equations satisfied by conformal blocks in genus one case [Zhu].
\par
It is not obvious to construct examples of VOA which satisfy $C_2$-finiteness condition.  Only few examples with $C_2$-finiteness condition are known.  One of them is a series of VOA called $W(p)$, $p=2, 3, \dots$, which was constructed by H. G. Kausch about twenty years ago [Kau].  It is very recently proved that VOA $W(p)$ satisfies $C_2$-finiteness conditions by D. Adamovic [AM2] and [CaN].  It is known that the abelian categories $W(p)$-mod are not semi-simple.
\par
Conformal field theory associated to a VOA $W(p)$ gives a logarithmic conformal field theory, because zero mode operator $T(0)$ of energy-momentum tensor is not diagonalizable, and therefore $N$-points functions may have logarithmic parts [Gab].
\par
Quite recently it is observed that $W(p)$ type conformal field theory appears as the scaling limit of some boundary conditions of the integrable lattice models, (c.f. Pearce et al. [PRZ], Bushlanov et al. [BFGT]).
\par
The purpose of this paper is to analyze the structure of the abelian category of $W(p)$-modules.  In order to solve this problem, we construct $W(p)$-modules $\mathcal{P}_s^\pm$, $s=1, \dots, p-1$, and prove that these are in fact projective $W(p)$-modules.
\par 
In the papers [FGST1], [FGST2], Feigin et al. conjectured that two abelian categories $W(p)$-mod and $\bar U_q(sl_2)$-mod are categorically equivalent as abelian categories.  By using the structure theorems of these projective modules $\mathcal{P}_s^\pm$ obtained in this paper we prove the conjecture of Feigin et al.
\par
The VOA $W(p)$ are constructed by using the free field realization of the Virasoro algebra with central charge $c_p=13-6(p+\displaystyle\frac{1}{p})$, $p=2, 3, \dots$, and screening operators.  There are two screening operators $Q_+(z)$ and $Q_-(z)$.  For each integer $1 \leqq s \leqq p-1$ and $\varepsilon=\pm$, we define the screening operator $Q_-^{[d_s^\varepsilon]}(z)$ from $Q_-(z)$ by using the iterated integral on a twisted local system.  The screening operators $Q_-^{[d_s^\varepsilon]}(z)$, $1 \leqq s \leqq p-1$ and $\varepsilon=\pm$, play a very important role in this paper.
\par
In \S2 we collect some structures of Fock space representations of the Virasoro algebra by using intertwing operators arising from $Q_+(z)$ and $Q_-^{[d_s^\varepsilon]}(z)$.  The results are well known in [FF1], [FF2], [Fel] and [TsK].
\par
Our VOA $W(p)$ are defined from the lattice vertex operator algebra $V_L$ using the screening operator $Q_-^{d_s^\varepsilon}(z)$.  The $C_2$-finiteness condition of $W(p)$ is already known in [Ada], [AM1] and [AM2].  These fact will be stated in \S3.
\par
In \S4 we construct $W(p)$-modules $\mathcal{P}_s^\pm$, $1\leqq s\leqq p-1$, by using the method of J.Fjelsted et al. [FFHST].  The $W(p)$-modules $\mathcal{P}_s^\pm$ ($1 \leqq s \leqq p-1$), which we will construct in this paper is obtained by deforming a $W(p)$-module $\mathcal{V}_s^+ \oplus \mathcal{V}_s^-$ by using screening operators $Q_-^{[d_s^\varepsilon]}(z)$.  The construction of $\mathcal{P}_s^\pm$, and an analysis of these $W(p)$-module are the most important parts of this paper.  The structure of $\mathcal{P}_s^\pm$, $1\leqq s\leqq p-1$, is described in Theorem 4-3 and Theorem 4-4.  These two theorems are a part of the main results of this paper.  By using these structure theorems, we determine completely Zhu's algebra $A_0(W(p))$ of the VOA $W(p)$, which is stated in Theorem 4-6.
\par
In \S5, we determine the Ext$^1$ group between simple objects $\mathcal{X}_s^\pm$, $s=1, \dots , p$, the results is given in Theorem 5-1 and 5-2.  By using the both theorems, we show that the abelian category $W(p)$-mod of $W(p)$-modules has the block decomposition
$$W(p)\mbox{-mod}=\bigoplus_{s=0}^p C_s.$$
The subcategories $C_0$ and $C_p$ are semi-simple consisting of simple objects $\mathcal X_0$ and $\mathcal X_p$, respectively.  But for $1\leqq s\leqq p-1$, $C_s$ is not semi-simple.  The set of simple objects of $C_s$ consists of two elements $\{ \mathcal{X}_s^+, \mathcal{X}_s^- \}$.  These results are stated at the first part of $\S 5$.
\par
Finaly we show that the $W(p)$ module $\mathcal{P}_s^\pm$, $1\leqq s\leqq p-1$, is a projective cover of $\mathcal{X}_s^\pm$, self-dual, and therefore injective.  On the module $\mathcal{P}_s^\pm$, the zero mode operator $T(0)$ of the energy-momentum tensor is not diagonizable.   To prove the projectivity of $\mathcal{P}_s^\pm$, we must know the detailed structure of $\mathcal{P}_s^\pm$, and show that Ext$^1$ groups between $\mathcal{P}_s^\pm$ and simple modules $\mathcal{X}_s^\pm$ are zero.
\par
On the very final step for $1\leqq s\leqq p-1$ we compute the endmorphism algebra,
$$B_s=\mbox{End}_{C_s}(\mathcal{P}_s), \qquad \mathcal{P}_s=\mathcal{P}_s^+ \oplus \mathcal{P}_s^ -.$$
The structure of $B_s$ is given in Theorem 6-4.  They are eight dimensional basic Artinian algebras, mutually  isomorphic to the basic algebra arising from $\bar U_q(sl_2)$-mod computed by Feigin [FGST1], [FGST2].
\par
The structures of these basic Artinian algebras are explicitly described.  This is stated in Theorem 6-2.
\par
Using the fact that two basic algebras coming from $W(p)$ and $\bar{U}_q(sl_2)$ are isomorphic, it is easy to prove by the conjectures of Feigin [FGST1], [FGST2]:
$$W(p)\mbox{-mod} \simeq \bar U_q(sl_2)\mbox{-mod}.$$
Since the abelian category $W(p)$-mod is not semi-simple, it is very interesting and important to analyze the structures of: (1) the fusion tensor products, (2) the monodromy representations mapping class group, bradi group, and (3) genus one and higher genus conformal blocks occuring in the conformal field theory associated with the VOA $W(p)$.  Having the results obtained in this paper we are now ready to study these problems.
\par
The second named author appreciates to H. Yamauchi for the discussion in the early stage of preparing of this paper.

\vspace{0.5in}
\section{Free Fields and Screening Charge Operators}
In this section we give a free field realization of the Virasoro algebra and intertwining operators.

\vspace{0.2in}
\subsection{Notations}

We fix an integer $p\geqq 2$, and set $\alpha_+=\sqrt{2p}$, $\alpha_-=-\sqrt{2/p}$ and $\alpha_0=\alpha_+ + \alpha_-$.  Then we have $\alpha_+\cdot \alpha_-=-2$,  $\alpha_+=-p\alpha_-$, $\displaystyle\frac{1}{\alpha_+}=-\displaystyle\frac{\alpha_-}{2}$, $\displaystyle\frac{\alpha_0^2}{2}=\displaystyle\frac{(p-1)^2}{p}$ and $\alpha_+\cdot \alpha_+=2p$.\\

Let us introduce an even integral lattice and its dual;
\begin{equation}\tag{2.1}
L=\mathbb Z\alpha_+,
\end{equation}
\begin{equation}\tag{2.2}
L^{\vee}=\mbox{Hom}_\mathbb Z(L, \mathbb Z)=\mathbb Z\cdot \frac{\alpha_-}{2}.
\end{equation}
For any integers $s$ and $n \in \mathbb Z$, we set
\begin{equation}\tag{2.3}
\lambda_s(n)=\frac{1-s}{2}\alpha_-+n\alpha_+,
\end{equation}
and
\begin{equation}\tag{2.4}
\Lambda_s=\{\lambda_s(n) ; n\in \mathbb Z\}.
\end{equation}
Then we see that, if $s_1-s_2\neq 0$ (mod $2p$), then $\Lambda_{s_1}\cap \Lambda_{s_2}=\emptyset$, $\lambda_{s+2p}(n)=\lambda_s(n+1)$, $\Lambda_{s+2p}=\Lambda_s$.  We have the $L$-orbit decomposition of $L^\vee$ as follows,
\begin{equation}\tag{2.5}
L^{\vee}=\bigsqcup_{-(p-1)\leqq s\leqq p}\Lambda_s.
\end{equation}
We set for $1\leqq s \leqq p-1$
\begin{align}\tag{2.6}
& \Lambda_s^+=\Lambda_s, \ \Lambda_s^-=\Lambda_{-s}, \\
& \Lambda_p^-=\Lambda_0, \ \Lambda_s^+=\Lambda_p. \notag
\end{align}
For each $\mu \in \mathbb C$ we set
\begin{equation}\tag{2.7}
h_\mu=\frac{1}{2} \left( \mu-\frac{1}{2}\alpha_0 \right)^2 -\frac{1}{8}\alpha_0^2.
\end{equation}
Then we have
\begin{equation}\tag{2.8}
h_\lambda=h_{\alpha_0-\lambda}=h_{\lambda^+},
\end{equation}
where for $\lambda \in \mathbb C$ we denote
\begin{equation}\tag{2.9}
\lambda^+=\alpha_0-\lambda.
\end{equation}
Then we have 
\begin{equation}\tag{2.10}
h_{\lambda_s(n)}=\frac{1}{4p} \left\{ (2np+s-p)^2-(p-1)^2 \right\},
\end{equation}
and the following formulas hold, for all $s$ ($0\leqq s\leqq p-1$);
$$h_{\lambda_{-s}}(n)=h_{\lambda_s(1-n)}, \quad \lambda_{-s}(n)+\lambda_s(1-s)=\alpha_0,$$
$$h_{\lambda_p}(n)=h_{\lambda_p(-n)}, \quad \lambda_p(n)+\lambda_p(-n)=\alpha_0.$$
We introduce the following sequence of numbers, for $1\leqq s\leqq p-1$, $n\geqq 0$;
\begin{align}\tag{2.11}
& h_s(2n)=h_{\lambda_s(-n)}=h_{\lambda_{-s}(n+1)}  = \frac{1}{4p}\{ ((2n+1)p-s)^2 - (p-1)^2 \}, \\
& h_s(2n+1)=h_{\lambda_s(n+1)}=h_{\lambda_{-s}(-n)} = \frac{1}{4p}\{ ((2n+1)p+s)^2 - (p-1)^2 \}, \notag \\
& h_0(n)=h_{\lambda_0(n+1)}=h_{\lambda_0(-n)}= \frac{1}{4p}\{( (2n+1)p)^2 - (p-1)^2 \}, \notag \\
& h_p(n)=h_{\lambda_p(n)}=h_{\lambda_p(-n)} = \frac{1}{4p}\{ (2np)^2- (p-1)^2 \}. \notag
\end{align}
Then we have an increasing and a decreasing series of rational numbers, for $0 \leqq s \leqq p$;
\begin{equation}\tag{2.12}
h_s(0)<h_s(1)<h_s(2)<\dots, \ h_s(n+1)-h_s(n)\in \mathbb Z_{\geqq 1}, \ n \geqq 0,
\end{equation}
\begin{equation}\tag{2.13}
h_{p-1}(1)>h_{p-2}(1)>\dots >h_1(1)>h_0(0)>h_1(0)>\dots >h_p(0).
\end{equation}
We see $h_p(0)= -\displaystyle\frac{(p-1)^2}{4p}$, $h_1(0)=0$ and $h_0(0)=\displaystyle\frac{p^2-(p-1)^2}{4p}$.
\par
We also define the sets of rational numbers, for $1\leqq s\leqq p-1$:
$$H_s^+=\{ h_s(2n) ; n\geqq 0 \}, H_s^-=\{ h_s(2n+1) ; n\geqq 0 \}, H_s=H_s^+\cup H_s^-,$$
$$H_p=H_p^+=\{ h_p(n) ; n\geqq 0 \},$$
$$H_0=H_p^-=\{ h_0(n) ; n\geqq 0 \}.$$
Then we see $H_s \cap H_{s'} = \emptyset$ if $s \ne s'$.  We set
\begin{equation}\tag{2.14}
H=\bigsqcup_{s=0}^p H_s.
\end{equation}

\vspace{0.5in}
\subsection{Free field realization of the Virasoro algebra}
First we introduce the free Bosonic field as follows:
\begin{equation}\tag{2.15}
\varphi(z)=\hat{a}+a(0)\log z-\sum_{n\ne 0}\frac{a(n)}{n}z^{-n},
\end{equation}
\begin{equation}\tag{2.16}
a(z)=\partial \varphi(z)=\sum_{n \in z}a(n)z^{-n-1}.
\end{equation}
This field is characterized by the operator product expansions (OPE)
$$\varphi(z) \varphi(w)\sim\log(z-w),$$
$$\partial \varphi(z)\partial \varphi(w)\sim \frac{1}{(z-w)^2}.$$
The operators $\hat{a}$ and $a(n)$ satisfy the following commutator relations
\begin{equation}\tag{2.17}
[a(n), a(n)]=m\delta_{m+n, 0} \ \mbox{id},
\end{equation}
$$[a(n), \hat{a}]=\delta_{n, 0}\hat{a}.$$
Set
$$\varphi_\pm(z)=\mp\sum_{n \geqq 1}\frac{a(\pm n)}{n}z^{\mp n},$$
then we have
$$\varphi(z)=\varphi_-(z)+\hat{a}+a(0)\log z+\varphi_-(z).$$
For each $\lambda \in \mathbb C$ we define the left and the right Fock module by following relations
\begin{align}\tag{2.18}
& F_\lambda \ni | \lambda\rangle \ne 0, \quad a(n)|\lambda\rangle=\delta_{n, 0}\lambda | \lambda\rangle, \quad n\geqq 0, \\
& F_\lambda^{\dagger}\ni \langle \lambda | \ne 0, \quad \langle \lambda | a(-n)=\delta_{n , 0}\lambda\langle \lambda|, \quad n\geqq 0. \notag
\end{align}
Then we have a unique non-degenerate pairing
\begin{equation}\tag{2.19}
\langle \ \ | \ \ \rangle:F_\lambda^\dagger \times F_\lambda \longrightarrow \mathbb C
\end{equation}
such that
$$\langle \lambda|\lambda\rangle=1, \quad \langle va(n) | u\rangle=\langle v | a(n)u \rangle,$$
for $n \in \mathbb Z, \ u \in F_\lambda, \ v \in F_\lambda^+$.
\par
Define the energy-momentum tensor
\begin{equation}\tag{2.20}
T(z)=\frac{1}{2}:\partial \varphi(z)^2:+\frac{\alpha_0}{2}\partial^2\varphi(z)=\sum_{n \in \mathbb Z}T(n)z^{-n-2},
\end{equation}
then we have OPE of the Virasoro field with central charge $c_p$;
\begin{equation}\tag{2.21}
T(z)T(w)\sim \frac{\frac{1}{2}c_p}{(z-w)^4}+\frac{2}{(z-w)^2}T(w)+\frac{1}{(z-w)}\partial_wT(w)
\end{equation}
where $c_p=13-6 \left(p+\displaystyle\frac{1}{p}\right)$ as usual.
\par
For each $\mu \in \mathbb C$ we set
\begin{equation}\tag{2.22}
V_\mu(z)=\  : e^{\mu\varphi(z)} : \ =e^{\mu\varphi_+(z)}e^{\mu\varphi_-(z)}e^{\mu\hat{a}}e^{\mu a(a)\log z}.
\end{equation}
Then we have
$$V_\mu(z):F_\lambda \longrightarrow F_{\lambda+\mu},$$
$$V_{\mu_1}(z_1)V_{\mu_2}(z_2)=(z_1-z_2)^{\mu_1\cdot \mu_2} : V_{\mu_1}(z_1)V_{\mu_2}(z_2),$$
and the following operator product expansion
\begin{equation}\tag{2.23}
T(z)V_\mu(w)\sim\frac{h_\mu}{(z-w)^2}V_\mu(w)+\frac{1}{(z-w)}\partial_wV_\mu(w).
\end{equation}
Here $h_\mu$ is defined in (2.7) and called conformal dimension of the field operator $V_\mu (z)$.

For $1\leqq s \leqq p-1$ we set
$$\mathcal{V}_s^\pm=\sum_{\lambda\in\Lambda_s^\mp}F_\lambda,$$
\begin{equation}\tag{2.24}
\mathcal{V}_p=\mathcal{V}_p^+=\sum_{\lambda\in\Lambda_p}F_\lambda, \quad \mathcal{V}_0=\mathcal{V}_p^-=\sum_{\lambda\in\Lambda_p^-}F_\lambda.
\end{equation}
Then for each $\lambda \in L$ we get
\begin{equation}\tag{2.25}
V_\lambda(z) \subset \mbox{End}_{\mathbb C}(\mathcal{V}_s^\pm)[[z,z^{-1}]]
\end{equation}
for all $1\leqq s \leqq p$.

\vspace{0.5in}
\subsection{Screening Operators}
Since $h_{\alpha_\pm}=1$, the field $Q_\pm(z)=V_{\alpha_\pm}(z)$ has the conformal dimension 1 with respect to the Virasoro field $T(z)$.  For each $\lambda \in L^\vee$, since $\alpha_+ \cdot \lambda \in \mathbb Z$ we see that
\begin{equation}\tag{2.26}
Q_+(z)\in \mbox{Hom} (F_\lambda, F_{\lambda+\alpha_+})[[z, z^{-1}]].
\end{equation}
Therefore
\begin{equation}\tag{2.27}
Q_+(0)=Q_+=\int dz \ Q_+(z) : F_\lambda \longrightarrow F_{\lambda+\alpha_+}
\end{equation}
commutes with $T(z)$.  Here we consider $\displaystyle\int dz$ as taking residue at $z=0$.
\par
While for $\lambda \in L^\vee$, $\alpha_- \cdot \lambda \not\in \mathbb Z$ in general, therefore $\displaystyle\int Q_-(z)dz : F_\lambda \rightarrow F_{\lambda+\alpha_-}$ cannot be defined for $\alpha_- \cdot \lambda \not\in \mathbb Z$.
\par
To construct an intertwining operator from $Q_-(z)$, we have to use an iterated integration on the twisted cycles.  To this end, we prepare some notations.
\par
For $d \geqq 1$, consider the product of screening operators
\begin{align}\tag{2.28}
Q_-(w_1) & \dots Q_-(w_d) \\
& = e^{d\hat{a}} \prod_{1 \leqq i<j \leqq d}(w_i-w_j)^{\frac{2}{p}} \prod_{j=1}^d w_j^{\alpha_-a(0)} e^{\alpha_-\sum_{j=1}^d \varphi_-(w_j)} e^{\alpha_-\sum_{j=1}^d \varphi_+(w_j)}, \notag
\end{align}
which is an acting on
$$\mathcal{V}_{[\lambda]} \longrightarrow \mathcal{V}_{[\lambda+d\alpha_-]}.$$
For $d \geqq 2$, define the complex manifold
\begin{equation}\tag{2.29}
X_d = \{ (w_1, \dots , w_d) \in (\mathbb C^\times)^d; \ w_i \ne w_j \},
\end{equation}
\begin{align}\tag{2.30}
Z_{d-1} =  \{ (\xi_1, \dots , \xi_{d-1}) \in ( \mathbb C \setminus \{0,1\} )^{d-1}; \ \xi_i \ne \xi_j \},
\end{align}
and define a map
\begin{align}\tag{2.31}
\mathbb C^\times \times Z_{d-1} & \rightarrow X_d\\
(w ; \xi_1, \dots , \xi_{d-1}) & \mapsto (w\xi_0, w\xi_1, \dots , w\xi_{d-1}) \notag
\end{align}
where we put $\xi_0=1$
\par
Then this map is $\mathbb C^\times$-equivariant isomorphism where
\begin{align}\tag{2.32}
& \lambda (w ; \xi_1, \dots , \xi_{d-1}) = (\lambda w ; \xi_1, \dots , \xi_{d-1}), \\
& \lambda (w_1, \dots , w_d) = (\lambda w_1, \dots , \lambda w_d).  \notag
\end{align}
\par
For each $\lambda \in L^\vee$ and $d \geqq 2$, we define multivalent functions respectively on $X_d$, $X_{d-1}$, by
\begin{equation}\tag{2.33}
\Phi_d^\lambda (w_1, \dots , w_d) = \prod_{1 \leqq i < j \leqq d} (w_i - w_j)^{\frac{2}{p}} \prod_{j=1}^d w_j^{\alpha_- \lambda}
\end{equation}
and
\begin{equation}\tag{2.34}
\bar \Phi_{d-1}^\lambda (\xi_1, \dots , \xi_{d-1}) = \prod_{0 \leqq i < j \leqq d-1} (\xi_i - \xi_j)^{\frac{2}{p}} \prod_{j=1}^{d-1} \xi_j^{\alpha_- \lambda}
\end{equation}
where we set $\xi_0=1$.  Then we have the formula
\begin{equation}\tag{2.35}
\Phi_d^\lambda (w_1, \dots , w_d) = \bar \Phi_{d-1}^\lambda (\xi_1, \dots , \xi_{d-1}) \cdot w^{\Delta_d(\lambda)}
\end{equation}
where
\begin{equation}\tag{2.36}
\Delta_d(\lambda)=\frac{1}{p}d(d-1)+d\alpha_- \lambda \in \frac{1}{p}\mathbb Z.
\end{equation}
\par
For $\lambda \in L^\vee$ and $d \geqq 2$, we denote $\bar S_{d-1}^{\lambda, *}$, the local system on $Z_{d-1}$ , determined by the monodromy of $\bar \Phi_{d-1}^{\lambda}$, and also denote $\bar S_{d-1}^{\lambda}$, the dual local system of $\bar S_{d-1}^{\lambda, *}$.  Then these local systems depend only on the class $[\lambda] \in L^\vee / L$ of $\lambda \in L^\vee$.  Therefore we can write $\bar S_{d-1}^{[\lambda] }$ etc.
\par
If we take an element $[\bar \Gamma] \in H_{d-1}(Z_{d-1}, \bar S_{d-1}^{[\lambda]})$, the integral
\begin{equation}\tag{2.37}
 \int_{[ \bar \Gamma]} Q_-(w\xi_0), \dots ,Q_-(w\xi_{d-1}) w^{d-1} d\xi_1 \dots d\xi_{d-1}
\end{equation}
define an element of
\begin{equation}\tag{2.38}
\mbox{Hom}_{\mathbb C} (\mathcal{V}_{[\lambda]} , \mathcal{V}_{[\lambda+d\alpha_{-1}]}) [ w, w^{-1} ],
\end{equation}
if $\Delta_d(\lambda) \in Z$.
\par
For $1 \leqq s \leqq p-1$ and $\varepsilon=\pm$, we define an integer $d_s^\pm$ by the following way;
$$d_s^+=p-s \ \mbox{and} \ d_s^-=s.$$
And we denote
$$\lambda_s^+=\lambda_{-s}(1) \in \Lambda_s^+, \quad \lambda_s^-=\lambda_s(0) \in \Lambda_s^-.$$
Then we have the following;
$$\Delta_{d_s^+}(\lambda_s^+) \in \mathbb Z, \quad \Delta_{d_s^-}(\lambda_s^-) \in \mathbb Z.$$
For $1 \leqq s \leqq p-1$ we define operators
$$Q_-^{[d_s^\pm]}(w) \in \mbox{Hom}_{\mathbb C} (\mathcal{V}_s^\pm , \mathcal{V}_s^\mp) [[ w, w^{-1} ]]$$
by the following way;
\begin{enumerate}
\item For $d_s^\pm=1$, we set
$$Q_-^{[d_s^\pm]}(z) = Q_-(z).$$
\item For $2 \leqq d_s^\pm \leqq p-1$, we set
$$Q_-^{[d_s^\pm]}(z) = \int_{\bar{\Gamma}} Q_-(w\xi_{d_s^\pm-1})w^{d_s^\pm-1} d\xi_1 \dots d\xi_{d_s^\pm-1}.$$
\end{enumerate}
We fix a cycle
$$[\bar{\Gamma}] \in H_{d_s^\pm-1}(Z_{d_s^\pm-1}, S_{d_s^\pm-1}^{[\lambda_s^\pm]})$$
which satisfies the following normalized conditions [TsK];
$$\int_{\bar{\Gamma}} \bar{\Phi}_{d_s^\pm-1}^{\lambda_s^\pm} (\xi_1 \dots \xi_{d_s^\pm-1}) d\xi_1 \dots d\xi_{d_s^\pm-1}=1.$$

\vspace{0.2in}
\noindent
{\bf Proposition 2-1} 
For $1 \leqq s \leqq p-1$ and $\varepsilon=\pm$, we have
\begin{enumerate}
\item $T(z)Q_-^{[d_s^\varepsilon]}(w)\sim \displaystyle\frac{1}{(z-w)^2}Q_-^{[d_s^\varepsilon]}(w)+\displaystyle\frac{1}{(z-w)^2}\partial_w Q_-^{[d_s^\varepsilon]}(w)$,
\item $[Q_-^{[d_s^\varepsilon]}(0), \ T(z)]=0$,
\item $[Q_+, \ Q_-^{[d_s^\varepsilon]}(0)]=0$.
\end{enumerate}

\begin{proof}
It can be proved in the standard way.  We give a proof of (3) only.
\par
Since $Q_+(z)Q_-(w)=(z-w)^{-2}:e^{\alpha_+\phi(z)+\alpha_-\phi(w)}:$,
we have 
$$[Q_+, Q_-(w)]=\frac{\alpha_+}{\alpha_++\alpha_-}\frac{\partial}{\partial w}V_{\alpha_+ +\alpha_-}(w).$$
So we get
\begin{align*}
\left[Q_+,Q_-^{[d_\lambda]}\right]
& = \left[Q_+,\int_\Gamma dw_1\dots dw_{d_\lambda}Q_-(w_1)\dots Q_-(wd_\lambda)\right] \\
& = \frac{1}{\alpha_+ + \alpha_-}\int_\Gamma d [\sum_{j=1}^{d_\lambda} (-1)^{j+1}V_{\alpha_-}(w_1)\dots V_{\alpha_+ + \alpha_-}(w_j) \dots V_{\alpha_-}(w_{d_\lambda}) \\
& \qquad \qquad \qquad \qquad \qquad \qquad \qquad dw_1\dots \overset{\vee}{dw_j} \dots dw_{d_\lambda} ] \\
& = 0 .
\end{align*}
\end{proof}

\vspace{0.5in}
\subsection{Abelian category $\mathcal{L}_{c_p}$-mod}
Let us consider the Virasoro algebra
\begin{equation}\tag{2.38}
\mathcal{L}=\sum_{n\in Z} \mathbb C T(n)\oplus \mathbb Cc
\end{equation}
with $c=c_p$ id.  Define Lie subalgebra as $\mathcal{L}_{>0}$ and $\mathcal{L}_{<0}$ of $\mathcal{L}$
\begin{equation}\tag{2.39}
\mathcal{L}_{>0}=\sum_{n \geqq 1} \mathbb C T(n), \quad \mathcal{L}_{<0}=\sum_{n \geqq 1} \mathbb C T(-n),
\end{equation}
and we define involtive anti-astromorphism of Lie algebra $\mathcal L$ by $\sigma(T(n))=T(-n)$ and $\sigma(c)=c$.
\par
Consider $\mathcal{L}$-module $M$ with the following properties.
\begin{enumerate}
\item $c=c_p$ \ id on $M$.
\item $M$ has the following decomposition $M=\displaystyle\sum_{h \in H(M)} M[h]$, where $H(M)=H_0(M)+\mathbb Z_{\geqq 0}$, for some finite subset $H_0(M)$ of $\mathbb C$, and for $h \in H(M)$, set
$M[h]=\{ m \in M : (T(0)-h)^nm=0$  for some  $n \geqq 0 \}$.  We further assume $\dim_{\mathbb C}M[h] < \infty $.
\end{enumerate}
\par
Let us denote
\begin{equation}\tag{2.40}
\mathcal{L}_{c_p}\mbox{-mod},
\end{equation}
the abelian category of left  $\mathcal{L}$-modules which satisfy the above conditions (1) and (2).
\par
For each $h \in \mathbb C$, let $M_h \ni | h \rangle$ be the Verma module of the highest weight $h$, the highest vector $| h \rangle$, and the central charge $c=c_p$ id, let $L_h$ be the irreducible quotient of $M_h$.
\par
These $\mathcal{L}_{c_p}$-module $M_h$, $L_h$, and Fock module $F_\lambda$ are objects of $\mathcal{L}_{c_p}$-mod.
\par
For each $\lambda \in \mathbb C$, there exists a unique $\mathcal{L}_{c_p}$-module map
\begin{equation}\tag{2.41}
M_{h_\lambda} \longrightarrow F_\lambda
\end{equation}
such that $| h_\lambda \rangle$ is mapped to $| \lambda \rangle$.
\par
The facts which we are going to use can be found in Feigin and Fuchs [FF1], [FF2], Felder [Fel], and Tsuchiya and Kanie [TsK].  By using Kac determinant formula for the Virasoro algebra, it is easy to show the following.

\vspace{0.2in}
\noindent
{\bf Proposition 2-2.} 
For $h\in \mathbb C\setminus H$, the Virasoro module $M_h$ is a simple object in $\mathcal{L}_{c_p}$-mod, where $H$ is defined in \S2-1, (2.10).

\vspace{0.2in}
\noindent
{\bf Proposition 2-3.}  
Fix $0 \leqq s \leqq p$, for $m, n \in \mathbb C$.  Then
\begin{enumerate}
\item 
\begin{equation}\notag
\mbox{Hom}_{C_s}(M_{h_s(m)}, M_{h_s(n)}) \simeq \left\{
\begin{gathered}
\mathbb C \quad m \geqq n \\
0 \quad m < n
\end{gathered} \right .
\end{equation}
\item For $m \geqq n$, the Virasoro sequence \\
$$0 \longrightarrow M_{h_s(m)} \longrightarrow M_{h_s(n)}$$
is exact.
\item For $n \geqq 0$, the Virasoro sequence \\
$$0 \longrightarrow M_{h_s(n+1)} \longrightarrow M_{h_s(n)} \longrightarrow L_{h_s(n)} \longrightarrow 0$$
is exact.
\end{enumerate}
\par
We define the following notations for the later use.  For $0 \leqq s \leqq p$ there exists a singular vector element
\begin{equation}\tag{2.42}
\eta_s | h_s(0)\rangle \in M_{h_s(0)} [h_s(1)]
\end{equation}
which is uniquely determined up to constant.  Where $\eta_s$ is an element
\begin{equation}\tag{2.43}
\eta_s \in U(\mathcal{L}_{<0})[s],
\end{equation}
we define
\begin{equation}\tag{2.44}
\eta_s^\vee = \sigma(\eta_s) \in U(\mathcal{L}_{>0})[-s].
\end{equation}

\vspace{0.2in}
\noindent
{\bf Proposition 2-4.}
\begin{enumerate}
\item For each $1 \leqq s \leqq p-1$ the followings hold:
  \begin{enumerate}
  \item $Q_+| \lambda_s(n)\rangle=0$, $Q_+^{2n} | \lambda_s(-n)\rangle\ne 0$ and $Q_+^{2n+1} | \lambda_s(-n)\rangle=0$ \ for \ $n \geqq 0$.
  \item $Q_+| \lambda_{-s}(n+1)\rangle=0$, $Q_+^{2n+1} | \lambda_{-s}(-n)\rangle\ne 0$ and $Q_+^{2n+2} | \lambda_{-s}(-n)\rangle=0$ \ for \ $n \geqq 0$. 
  \end{enumerate}
\item $Q_+| \lambda_0(n+1)\rangle=0$, $Q_+^{2n+1} | \lambda_0(-n)\rangle\ne 0$ and $Q_+^{2n+2} |\lambda_0(-n)\rangle= 0$ \ for \ $n \geqq 0$.
\item $Q_+| \lambda_p(n)\rangle=0$, $Q_+^{2n} | \lambda_p(-n)\rangle\ne0$ and $Q_+^{2n+1} | \lambda_p(-n)\rangle=0$ \ for \ $n \geqq 0$.
\end{enumerate}

\vspace{0.2in}
\noindent
{\bf Proposition 2-5.}  
For $1 \leqq s \leqq p-1$, we have:
\begin{enumerate}
\item $Q_-^{[s]} | \lambda_s(-n))=0$, $\quad Q_-^{[s]} | \lambda_s(n+1)\rangle\ne0$ for $n\geqq 0.$
\item $Q_-^{[p-s]} | \lambda_{-s}(-n)\rangle=0$, $\quad Q_-^{[p-s]} | \lambda_{-s}(n+1)\rangle\ne0$ for $n \geqq 0.$
\item $Q_-^{[p-s]} | \lambda_{-s}(1)\rangle=c|\lambda_s(0)\rangle \quad (c \ne 0).$
\end{enumerate}

\vspace{0.2in}
\noindent
{\bf Proposition 2-6.} 
 We have the following exact sequences of Virasoro modules with $c=c_p$.
\begin{enumerate}
\item For $1 \leqq s \leqq p-1$, $n \geqq 0$:
  \begin{enumerate}
  \item $0\longrightarrow M_{h_s(2n+1)}\longrightarrow M_{h_s(2n)}\longrightarrow F_{\lambda_s(-n)}$,
  \item $0\longrightarrow M_{h_s(2n+3)}\longrightarrow M_{h_s(2n+1)}\longrightarrow F_{\lambda_s(n+1)}$,
  \item $0\longrightarrow M_{h_s(2n+2)}\longrightarrow M_{h_s(2n)}\longrightarrow F_{\lambda_{-s}(n+1)}$,
  \item $0\longrightarrow M_{h_s(2n+2)}\longrightarrow M_{h_s(2n+1)}\longrightarrow F_{\lambda_{-s}(-n)}$.
  \end{enumerate}
\item For $s=0$, $n\geqq 0$:
  \begin{enumerate}
  \item $0\longrightarrow M_{h_0(n+1)}\longrightarrow M_{h_0(n)}\longrightarrow F_{\lambda_0(-n)}$,
  \item $0\longrightarrow M_{h_0(n+1)}\longrightarrow M_{h_0(n)}\longrightarrow F_{\lambda_0(n+1)}$.
  \end{enumerate}
\item For $s=p$, $n\geqq 0$:
  \begin{enumerate}
  \item $0\longrightarrow M_{h_p(n+1)}\longrightarrow M_{h_p(n)}\longrightarrow F_{\lambda_p(-n)}$,
  \item $0\longrightarrow M_{h_p(n+1)}\longrightarrow M_{h_p(n)}\longrightarrow F_{\lambda_p(n)}$.
  \end{enumerate}
\end{enumerate}

\vspace{0.2in}
As a consequence we obtain the so-called Felder complex [Fel].

\vspace{0.2in}
\noindent
{\bf Theorem 2-7.} 
 For $1\leqq s \leqq p-1$, the following is exact sequence of Virasoro modules
$$\overset{Q_-^{[s]}}{\longrightarrow} \mathcal{V}_s^-\overset{Q_-^{[p-s]}}{\longrightarrow} \mathcal{V}_s^+\overset{Q_-^{[s]}}{\longrightarrow} \mathcal{V}_s^- \longrightarrow \dots$$

\vspace{0.2in}
We define
$$\mathcal{X}_s^\pm=\ker Q^{[d_s^\mp]}, \quad Q^{[d_s^\mp]} : \mathcal{V}_s^\mp \longrightarrow \mathcal{V}_s^\pm.$$
Then we also have exact sequences of Virasoro modules
$$0 \longrightarrow \mathcal{X}_s^\mp \longrightarrow \mathcal{V}_s^\pm \longrightarrow \mathcal{X}_s^\pm \longrightarrow 0.$$
Virasoro modules $\mathcal{X}_s^\pm$ and $\mathcal{X}_p^\pm$ are decomposed into the sum of Virasoro submodules.

\vspace{0.2in}
\noindent
{\bf Theorem 2-8.}
\begin{enumerate}
\item For $1\leqq s \leqq p-1$,\\
$\mathcal{X}_s^+=\displaystyle\sum_{n=0}^\infty \sum_{m=0}^{2n}U(\mathcal{L})Q_+^m | \lambda_s(-n)\rangle$, \\

\qquad$U(\mathcal{L})Q_+^m|\lambda_s(-n)\rangle \simeq L_{h_s(2n)}$. 

\item For $1\leqq s\leqq p-1$,\\
$\mathcal{X}_s^-=\displaystyle\sum_{n=0}^\infty \sum_{m=0}^{2n+1}U(\mathcal{L})Q_+^m | \lambda_{-s}(-n)\rangle$, \\

\qquad$U(\mathcal{L})Q_+^m|\lambda_{-s}(-n)\rangle \simeq L_{h_s(2n+1)}$.

\item $\mathcal{X}_p^+=\displaystyle\sum_{n=0}^\infty \sum_{m=0}^{2n}U(\mathcal{L})Q_+^m | \lambda_p(-n)\rangle$, \\

\qquad$U(\mathcal{L})Q_+^m|\lambda_p(-n)\rangle \simeq L_{h_p(n)}$.

\item $\mathcal{X}_p^-=\displaystyle\sum_{n=0}^\infty \sum_{m=0}^{2n+1}U(\mathcal{L})Q_+^m | \lambda_0(-n)\rangle$, \

\qquad$U(\mathcal{L})Q_+^m|\lambda_0(-n)\rangle \simeq L_{h_0(n)}$.
\end{enumerate}

\vspace{0.2in}
\subsection{Block structure of $\mathcal{L}_{c_p}$-mod}

\vspace{0.2in}
Consider the decomposition of $\mathbb C=\bigsqcup_{b\in B} b$, where $b=\{h\}$ for $h\in \mathbb C\setminus H$ or $b=H_s$, $0 \leqq s \leqq p$.  Let us consider the abelian subcategory $C_b(\mathcal{L}_{c_p})$ of $\mathcal{L}_{c_p}$-mod, which is parametrised by $b \in B$ as follows.
\begin{enumerate}
\item $b=\{h\}$, $h\in \mathbb C\setminus H$, then $M \in \mathcal{L}_{c_p}$-mod belongs to $C_b(\mathcal{L}_{c_p})$ if and only if $M$ is a direct sum of the Verma modules $M_h$.
\item $b=H_s$, $0\leqq s \leqq p$, then $M \in \mathcal{L}_{c_p}$-mod belongs to $C_b(\mathcal{L}_{c_p})$ if and only if the irreducible sub-quotient of $M$ is isomorphic to $L_{h_s(n)}$ ($n \geqq 0$).
\end{enumerate}

\vspace{0.2in}
\noindent
{\bf Theorem 2-9.} 
The abelian category $\mathcal{L}_{c_p}$-mod has the following decomposition of abelian category
$$\mathcal{L}_{c_p}\mbox{-mod}=\bigoplus_{b \in B}C_b(\mathcal{L}_{c_p}).$$
For $b \ne b'$, $b, b' \in B$ and $M \in C_b(\mathcal L_{c_p})$, $N \in C_{b'}(\mathcal L_{c_p})$, we have the following facts,
$$\mbox{Ext}_{\mathcal L_{c_p}}^i(M, N) = 0. \qquad i=0,1,\dots .$$

\vspace{0.2in}
The homological properties of the abelian category $C_s(\mathcal{L}_{c_p})=C_{H_s}(\mathcal{L}_{c_p})$, $0 \leqq s\leqq p$, are very important in this paper.  So we would like to explain the required results.  We can not find in the literalness, but the results can be proved by using Kac determinant formula for the Virasoro algebra, and proofs of structure theorems of Virasoro modules are due to Feigin et al. [FF1], [FF2].
\par
At first we fix $s$, $1 \leqq s\leqq p-1$, and consider the abelian category $C_s(\mathcal{L}_{c_p})$.  We use the following notations
\begin{align}\tag{2.45}
& M_n=M_{h_s(n)}, \quad n \geqq 0, \\
& L_n=M_{h_s(n)} / M_{h_s(n+1)}, \quad n \geqq 0, \notag \\
& L_n^{(1)}=M_{h_s(n)} / M_{h_s(n+2)}, \quad n \geqq 0, \notag \\
& L_n^{(1)\vee}=D(L_n^{(1)}), \quad n \geqq 0. \notag
\end{align}

\vspace{0.2in}
\noindent
{\bf Proposition 2-10.}  
For each $s$, $1 \leqq s\leqq p-1$, we have the following.
\begin{enumerate}
\item The set of equivalense classes of simple objects in $C_s(\mathcal{L}_{c_p})$ are $\{ L_n : n \geqq 0 \}.$
\item The module $L_n$ are self dual $D(L_n) \simeq L_n$.
\item For $m, n \in \mathbb Z_{\geqq 0}$, we have
$$
\mbox{Ext}^1(L_m, L_n) \simeq
\begin{cases}
C & m=n \pm 1, \\
0 & \mbox{otherwise}.
\end{cases}
$$
\item $\mbox{Ext}^1(L_n, L_{n+1}) \ni [L_n^{(1)}] \ne 0$, \\
$\mbox{Ext}^1(L_{n+1}, L_n) \ni [L_n^{(1)\vee}] \ne 0$.
\end{enumerate}

\vspace{0.2in}
Now we restrict our attention to $1 \leqq s \leqq p-1$, and fix the following highest weight vectors;
\begin{align}\tag{2.46}
& u \in L_0[h_s(0)], \ u_0^{(1)} \in L_0^{(1)}[h_s(0)], \\
& v \in L_1[h_s(1)], \notag \\
& \eta_s(u_0^{(1)})=v \quad \mbox{in} \ L_0^{(1)}. \notag
\end{align}
\par
Then we have the following exact sequences;
$$0 \longrightarrow L_1 \longrightarrow L_0^{(1)} \longrightarrow L_0 \longrightarrow 0$$
$$\qquad \quad u_0^{(1)} \ \mapsto \ u$$

\vspace{0.2in}
\noindent
{\bf Proposition 2-11.}  
Fix $s$, $1 \leqq s\leqq p-1$, then the following hold.
\begin{enumerate}
\item $\mbox{Ext}^1(L_0, L_0^{(1)})=0$, \  $\mbox{Ext}^1(L_0^{(1)\vee}, L_0)=0$.
\item $\mbox{Ext}^1(L_0^{(1)}, L_0) \simeq \mathbb C$, \  $\mbox{Ext}^1(L_0, L_0^{(1)\vee}) \simeq \mathbb C$.
\item Fix a generator $K^{(1)}$ of $\mbox{Ext}^1(L_0^{(1)}, L_0) \simeq \mathbb C$.  Then the following hold.
\begin{enumerate}
\item $D(K^{(1)}) \simeq K^{(1)}$.
\item $K^{(1)}$ is a generator of $\mbox{Ext}^1(L_0, L_0^{(1)\vee}) \simeq \mathbb C$,
$$[K^{(1)}] \in \mbox{Ext}^1(L_0, L_0^{(1)\vee}).$$
\end{enumerate}
\item  We can take elements $u_0, u_1 \in K^{(1)}[h_s(0)]$ and $v_0 \in K^{(1)}[h_s(1)]$ with the following properties
\begin{align}\tag{2.47}
& v_0=\eta_s(u_0),  \\
& u_1=\eta_s^\vee(v_0) \in L_0 \subseteq K_0^{(1)}, \notag 
\end{align}
$$\ 0 \longrightarrow L_0 \longrightarrow K^{(1)} \longrightarrow L_0^{(1)} \longrightarrow 0,$$
$$\qquad \qquad u_0 \ \mapsto \ u_0^{(1)}$$
$$\qquad \quad v_0 \ \mapsto \ v$$
Then the following relations hold;
\begin{equation}\tag{2.48}
(T(0)-h_s(0))u_0=cu_1, \quad c \ne 0.
\end{equation}
\end{enumerate}

\vspace{0.2in}
\noindent
{\bf Proposition 2-12.}  
Fix $s$, $1 \leqq s\leqq p-1$, then the following hold.
\begin{enumerate}
\item $\mbox{Ext}^1(L_0^{(1)}, L_1)=0$, \  $\mbox{Ext}^1(L_1, L_0^{(1)\vee})=0$.
\item $\mbox{Ext}^1(L_1, L_0^{(1)}) \simeq \mathbb C$, \  $\mbox{Ext}^1(L_0^{(1)\vee}, L_1) \simeq \mathbb C$.
\item Fix the generator $K_{(1)}$ of $\mbox{Ext}^1(L_1, L_0^{(1)}) \simeq \mathbb C$, the following facts hold.
\begin{enumerate}
\item $D(K_{(1)}) \simeq K_{(1)}$.
\item $K_{(1)}$ is a generator of $\mbox{Ext}^1(L_0^{(1)\vee}, L_1) \simeq \mathbb C$,
$$[K_{(1)}] \in \mbox{Ext}^1(L_0^{(1)\vee}, L_1).$$
\end{enumerate}
\item  We can take elements $u_0$, $v_0$, and $v_1 \in K_{(1)}$ \\
such that
$$u_0 \in L_0^{(1)}[h_s(0)]=K_{(1)}[h_s(0)],$$
$$v_1 \in L_0^{(1)}[h_s(1)] \subseteq K_{(1)}[h_s(1)], \ v_0 \in K_{(1)}[h_s(1)].$$
$$\eta_s^\vee(v_0)=u_0, \ \eta_s(u_0)=v_1$$
$$0 \longrightarrow L_0^{(1)} \longrightarrow K_{(1)} \longrightarrow L_1 \longrightarrow 0,$$
$$\qquad \qquad v_0 \mapsto v.$$
Then the following relations hold.
$$(T(0)-h_s(1))v_0=cv_1, \quad c \ne 0.$$
\end{enumerate}

\vspace{0.2in}
Finally we obtain the following proposition.

\vspace{0.2in}
\noindent
{\bf Proposition 2-13.}  
The abelian category $C_{H_0}(\mathcal L_{c_p})$ and $C_{H_p}(\mathcal L_{c_p})$ are semi-simple.

\vspace{0.5in}
\section{The Triplet VOA $W(p)$}
In this section, we define the so-called triplet VOA $W(p)$ and show that it satisfies Zhu's $C_2$-finiteness condition [AM2].

\vspace{0.5in}
\subsection{Vertex Operator algebras} 
In this paper the notion of vertex operator algebra (VOA) plays an important role.  For definitions and properties of VOA, we follow [FrB] and [Kac].  We use the notations of [NaT].
\par
Roughly speeking, a vertex operator algebra is a quadruple $(V, |0\rangle, T, J)$ such that
\begin{equation}\tag{3.1}
V=\bigoplus_{\Delta\geqq 0}V[\Delta]
\end{equation}
which is a $\mathbb Z_{\geqq 0}$ graded $\mathbb C$-vector space with the properties
$V[0]=\mathbb C | 0\rangle \ne 0$, $\dim_{\mathbb C}V[\Delta]<\infty$, and with an distinguished element $T \in V[2]$, $T \ne 0$.  The element $| 0 \rangle$ is called the vacum element and the element $T$ is called the Virasoro element.
\par
There exists a degree preserving linear map
\begin{equation}\tag{3.2}
J : V \longrightarrow \mbox{End}(V)[[z,z^{-1}]],
\end{equation}
$$A \mapsto J(A,z),$$
where we set degree of $z=-1$.  These must satisfy some compatibility conditions.  The most important properties are the locality of any two operators $J(A,z)$ and $J(B,w)$, and their operator product expansions (OPE).  For details of OPE, we refer [MaN].
\par
For each $A\in V[\Delta]$, we denote
\begin{equation}\tag{3.3}
J(A , z)=\sum_{n \in \mathbb Z} A_{(n)}z^{-n-1}=\sum_{n \in \mathbb Z}  A[n]z^{-n-\Delta}
\end{equation}
where $A_{(n)}=A[n-\Delta+1]$, $A[n]=A_{(n+\Delta-1)}$, $\deg A_{(n)}=-n+\Delta-1$ and $\deg A[n]=-n$.   Sometimes we write $J_n(A)=A[n]$.
\par
We denote
\begin{equation}\tag{3.4}
J(T,z)=T(z)=\sum_{n \in \mathbb Z}T(n)z^{-n-2}
\end{equation}
Then $\deg T(n)=-n$, and we have the following operator product expansion (OPE);
\begin{equation}\tag{3.5}
T(z)T(w)\sim \frac{\frac{1}{2}c}{(z-w)^4}+\frac{2}{(z-w)^2}T(a)+\frac{1}{(z-w)}\partial_wT(w)
\end{equation}
where $c$ is a some complex number.  The operator $T(z)$ is called the energy-momentum tensor.
\par
For $A, B \in V$, we denote the OPE for $J(A,z)$ and $J(B,z)$ by the following way.
\begin{align}\tag{3.6}
J(A,z)J(B,z)
& = \sum_{n \in \mathbb Z}J(J_n(A)B,w)(z-w)^{-n-\Delta_A} \\
& = \sum_{n \in \mathbb Z}J(A_{(n)}B,w)(z-w)^{-n-1}. \notag
\end{align}
And sometimes we use the following notations
\begin{equation}\tag{3.7}
J(A_{(n)}B,w)=(A_{(n)}B)(w),
\end{equation}
$$(A_{(n)}B)(w)=\mbox{Res}_{z=w} (z-w)^n A(z)B(w) \ dz.$$
\par
The representation $(M,J^M)$ of a VOA is a degree preserving linear map
\begin{equation}\tag{3.8}
J^M : V \longrightarrow \mbox{End}(M)[[z,z^{-1}]]
\end{equation}
such that $M=\sum_{h\in H(M)} M[h]$ with $M[h]=\{ m \in M ; (T(0)-h)^n m=0$ for some $n\geqq 0\}$ and some compatibility conditions.  In this paper we assume that $H(M) = H_0(M)+\mathbb Z_{\geqq 0}$ for a finite set $H_0(M)$, and also assume that for any $h \in H(M)$, $\dim_{\mathbb C}M[h]<\infty$.  In general $\dim_{\mathbb C}M[h]<\infty$ is a too strong condition.  However, since in this paper we mainly deal with VOA's which satisfy $C_2$-finiteness conditions, this condition is not restrictive.  We denote
\begin{equation}\tag{3.9}
V\mbox{-mod}
\end{equation}
the abelian category of left $V$-modules which satisfy the above conditions.
\par
For a VOA $V$, its universal enveloping algebra
\begin{equation}\tag{3.10}
U(V)=\sum_d U(V)[d]
\end{equation}
is introduced in [FrZ], [NaT] and [MNT].
\par
The algebra $U(V)$ is a degreewise completed linear topological algebra generated by $A[n]$, $A\in V$ and $n\in \mathbb Z$, $\deg A[n]=-n$.  A representation of VOA $V$ is a representation of $U(V)$, and vice versa.
\par
We define a subalgebra $F_0(U(V))=\sum_{d \leqq 0}U(V)[d]$ of $U(V)$ and a closed left ideal $I_0(V)$ of $U(V)$, which is generated by $\sum_{d \leqq -1}U(V)[d]$.
Then $F_0(U(V))\cap I_0(V)$ is a closed two-sided ideal of $F_0(U(V))$.  The Zhu's algebra $A_0(V)$ of $V$ is defined as the quotient algebra of $F_0(U(V))$,
\begin{equation}\tag{3.11}
A_0(V)=F_0(U(V)) / F_0(U(V))\cap I_0(V).
\end{equation}
For any $A \in V$, let $[A[0]]$ be the element of $A_0(V)$ represented by $A[0]$ mod $I_0(V)$.  The algebra $A_0(V)$ also can be defined as a quotient space of $V$ itself [Zhu].
\par
The algebra $A_0(V)$ is called zero mode algebra or Zhu's algebra of VOA $V$.
\par
For each $M \in V$-mod, define
\begin{equation}\tag{3.12}
HW(M)=\{ m \in M : J_n(A)m=0, \quad n \geqq 1, A \in V\}.
\end{equation}
Then the Zhu algebra $A_0(V)$ act on $HW(M)$.
\par
Here we introduce one important notion called Zhu's $C_2$-finiteness condition.
\par
For each vertex operator algebra $V$, we define a graded subspace $C_2(V)$ of $V$ by
\begin{equation}\tag{3.13}
C_2(V)=\mbox{ span of }\{ A_{(n)}B : A, B \in V, \ n \leqq -2 \}.
\end{equation}
Then the quotient space
\begin{equation}\tag{3.14}
\mathfrak{p}(V)=V/C_2(V)
\end{equation}
is graded, and has a Poisson algebra structure defined by for any $A, B \in V$;
\begin{align}\tag{3.15}
& [A]\cdot [B]=[A_{(-1)}B] , \\
& \{ [A], [B] \}=[A_{(0)}B] \notag
\end{align}
where $[A]$ denote the equivalent class of $A$ in $\mathfrak{p}(V)$.

\vspace{0.2in}
\noindent
{\bf Difinition 3-1.}  
The following condition of $V$ is called Zhu's $C_2$-finiteness condition;
\begin{equation}\tag{3.16}
\dim_{\mathbb C}\mathfrak{p}(V)<\infty .
\end{equation}

\vspace{0.2in}
When $\dim_{\mathbb C}\mathfrak{p}(V)<\infty$, we are able to prove $\dim_{\mathbb C}A_0(V)\leqq \dim_{\mathbb C}\mathfrak{p}(V)$.  Thus in this case we would like to study $A_0(V)$-mod that is the abelian category consisting of all finite dimensional left $A_0(V)$-modules.  Then the covariant functor $HW$ maps any $V$-modules to a $A_0(V)$-module,
\begin{equation}\tag{3.17}
HW:V\mbox{-mod} \longrightarrow A_0(V)\mbox{-mod},
\end{equation}
and it has the adjoint functor
\begin{equation}\tag{3.18}
X \in A_0(V)\mbox{-mod} \mapsto U(V)\underset{F_0(U(V))}{\otimes}X \in V\mbox{-mod},
\end{equation}
where $X$ is considered as a $F_0(U(V))$-module through the map $F_0(U(V)) \rightarrow A_0(V)$.
\par
The following important theorem is due to [FrZ], [MNT].

\vspace{0.2in}
\noindent
{\bf Theorem 3-2.}  
Suppose that $V$ satisfies Zhu's $C_2$-finiteness condition.  Then we have:
\begin{enumerate}
\item The abelian category $V$-mod is Artinian and Noetherian.
\item The number of equivalence classes of simple $V$-modules is finite.
\item The number of simple $A_0(V)$-modules is equal to the number of simple $V$-modules.
\end{enumerate}

\vspace{0.5in}
\subsection{The lattice vertex operator algebra $V_L$}
Define
\begin{equation}\tag{3.19}
V_L=\sum_{\lambda \in L}F_\lambda
\end{equation}
and set $T=\displaystyle\frac{1}{2}a(-1)^2|0\rangle-\frac{\alpha_0}{2}a(-2)|0\rangle \in V_L$, then the following is well known [FrB].

\vspace{0.2in}
\noindent
{\bf Theorem 3-3.}
\begin{enumerate}
\item There exists a unique vertex operator algebra structure on $V_L$ such that
\begin{align*}
& J(a(-1) | 0\rangle : z)=a(z), \\
& J(|\lambda\rangle:z)=V_\lambda(z)\quad \mbox{for} \ \lambda \in V_L,\\
& J(T:z)=T(z).
\end{align*}
\item For each $1\leqq s\leqq p$, $\mathcal{V}_s^\pm$ is an irreducible $V_L$-module.
\item The abelian category of $V_L$-modules is semi-simple and its inequivalent simple objects are $\mathcal{V}_s^\pm$, $1\leqq s\leqq p$.
\end{enumerate}

\vspace{0.2in}
Then $(F_0, |0\rangle, T, J)$ is a vertex operator subalgebra of $V_L$.  We remark that Fock space $F_0$ is not the one which appears in the filtration of $U(V)$ and think this may not make any confusions.  The VOA $F_0$ is generated by fields $a(z)$ and the associated Virasoro field is $T(z)=\displaystyle\frac{1}{2}:a(z)^2:+\alpha_0 / 2 \ \partial a(z)$.
\par
Note that $V(\mathcal{L}_{c_p}): = U(\mathcal{L}) | 0 \rangle \subseteq F_0$ contains $ | 0 \rangle$ and $T$, therefore $V(\mathcal{L}_{c_p})$ is a sub VOA of $F_0$.  The abelian category $V(\mathcal{L}_{c_p})$-mod is nothing but the abelian category $\mathcal{L}_{c_p}$-mod.  The $\mathcal{L}$-module $V(\mathcal{L}_{c_p})$ is isomorphic to $L_{h_1(0)}$ as $\mathcal{L}$-modules.  Note that $h_1(0)=0$.

\vspace{0.5in}
\subsection{Duality in $V$-mod}
The duality functor in VOA was introduced in [FHL].  The universal enveloping algebra $U(V)$ has an involutive anti-automorphism of the topological algebra $U(V)$:
\begin{equation}\tag{3.20}
\sigma:U(V) \longrightarrow U(V), \quad (\mbox{white} \ \sigma(A)=A^\sigma \ \mbox{for short}),
\end{equation}
such that $\sigma(U(V)[d])=U(V)[-d]$.  For $A \in V[\Delta_A]$, we define
\begin{equation}\tag{3.21}
J^\sigma(A;z)=\sum_n A[n]^\sigma z^{-n-\Delta_A},
\end{equation}
which is given by
\begin{equation}\tag{3.22}
J^\sigma(A ; z)=J(e^{zT(1)}(-z^{-2})^{T(0)}A; z^{-1})
\end{equation}
For $M \in V$-mod, its dual $D(M) \in V$-mod is defined by
\begin{equation}\tag{3.23}
D(M)=\sum_{h \in H(M)}\mbox{Hom}_{\mathbb C}(M[h], \mathbb C)
\end{equation}
as $\mathbb C$-vector space and the action is defined by
\begin{equation}\tag{3.24}
\langle A \phi, u \rangle=\langle \phi, A^\sigma u \rangle
\end{equation}
for all $A \in U(V)$, $\phi \in D(M)$ and $u \in M$.
\par
The following gives the duality of the VOA $V_L$.

\vspace{0.2in}
\noindent
{\bf Proposition 3-4.}  We have
\begin{enumerate}
\item $a(n)^\sigma=-a(-n)+\delta_{n, 0}\alpha_0$ id, $V_{\lambda}(n)^\sigma=-V_{\lambda}(-n)$, \ $\lambda \in L $ \\
$\sigma(T(n))=T(-n)$ \ for $n\in \mathbb Z$.
\item $D(\mathcal{V}_s^\pm)=\mathcal{V}_s^\mp$ \ for $1 \leqq s \leqq p-1$.
\item The sub-VOA $F_0$ is closed under the duality and $D(F_\lambda)=F_{\alpha_0-\lambda}$ for any $\lambda\in \mathbb C$.
\end{enumerate}

\vspace{0.2in}
\subsection{Construction of $W(p)$}
Recall that $\mathcal{V}_1^-=V_L$ carries a VOA structure.  Then the intertwining operator $Q_-^{[1]}=Q_-$ defines a subspace
$$W(p)=\mathcal{X}_1^+=\ker \ (Q_- : V_L \longrightarrow \mathcal{V}_1^+).$$
The space $W(p)$ contains $| 0 \rangle$ and $T$.  Thus $W(p)=(W(p), | 0 \rangle, T, J)$ defines a sub VOA of $V_L$.  This VOA W(p) is called the triplet VOA.
\par
We denote
$$W^-=|-\alpha_+\rangle, \quad W^0=Q_+|-\alpha_+\rangle, \quad W^+=Q_+^2|-\alpha_+\rangle.$$
Then we see that $W^a \in W(p)$, $a \in \{\pm, 0\}$ by Proposition 2-7 and that its conformal dimension is $2p-1$.

\vspace{0.2in}
\noindent
{\bf Proposition 3-5.}  
The VOA $W(p)$ is generated by $T(z)$, $W^a(z)=J(W^a:z)$, $a \in \{\pm, 0\}$ as a VOA.

\vspace{0.2in}
For the proof we refer the readers to [AM2].

\vspace{0.2in}
\subsection{$C_2$-finiteness of $W(p)$}
We define $W_0(p)=W(p)\cap F_0$.  Then $W_0(p)$ is a sub-VOA both of $F_0$ and $W(p)$.  It is easy to see that $T, W^0 \in W_0(p)$.

\vspace{0.2in}
\noindent
{\bf Proposition 3-6.}   
The VOA $W_0(p)$ is generated by $T(z)$ and $W^0(z)$ as a VOA.

\vspace{0.2in}
For the proof we see [Ada].

\vspace{0.2in}
Now we denote Zhu's algebra of $W_0(p)$ by $A_0(W_0(p))$.  Then by Proposition 3-6 $A_0(W_0(p))$ is a quotient algebra of polynomial ring $C[[T(0)], [W^0(0)]]$.
\par
Then the following important proposition is proved in [Ada].

\vspace{0.2in}
\noindent
{\bf Proposition 3-7. } 
Zhu's algebra $A_0(W_0(p))$ is isomorphic to $\mathbb C[[T(0)],[W^0(0)]] / \langle G\rangle$ where $\langle G\rangle$ is the ideal generated by an element
$$G=([W^0(0)]^2-c([T(0)]-h_p(0))\prod_{s=1}^{p-1}([T(0)]-h_s(0))^2$$
where $c=\displaystyle\frac{(4p)^{2p-1}}{((2p-1)!)^2}$.

\vspace{0.2in}
Now recall that
$$\mathfrak{p}=\mathfrak{p}(W(p))=W(p) / C_2W(p)$$
where
$$C_2(W(p)) \equiv \{ A_{(n)}B : A, B \in W(p) \quad n \leqq -2\}.$$

It is known that the associative algebra $A_0(W(p))$ has a filtration $G_{\bullet} A_0(W(p))$, so that we have a surjection
$$\mathfrak{p}(W(p)) \longrightarrow Gr_\bullet^GA_0(W(p)) \longrightarrow 0$$
as Poisson algebras, by $[T(0)] \mapsto [T]$, $[W^a(0)] \mapsto [W^a]$, $a\in \{\pm, b \}$ [MNT].
\par
The following proposition and corollary are proved in [AM2].

\vspace{0.2in}
\noindent
{\bf Proposition 3-8.}  
There exist the following relations on the Poisson algebra $\mathfrak{p}(W(p))$.
\begin{enumerate}
\item $[W^\pm]^2=0$, \quad $[W^0]^2+[W^-][W^+]=0$,\quad$[W^0][W^\pm]=0$.
\item $[W^0]^2=c[T]^{2p-1}$ \ ($c \ne 0$).
\item $[T]^p[W^a]=0$ \ ($a\in \{\pm, 0\}$).
\item $[T]^{3p-1}=0$.
\item $\{ [T], [W^a] \} =c_a[T]^p \ (c_0\ne 0, c_\pm=0).$
\item Other Poisson brackets are zero.
\end{enumerate}

\vspace{0.2in}
\noindent
{\bf Corollary 3-9.} 
 We have the following:
\begin{enumerate}
\item $\dim \mathfrak{p}(W(p)) \leqq 6p-1.$
\item $\dim A_0(W(p)) \leqq 6p-1.$
\item $W(p)$ satisfies Zhu's $C_2$-finiteness condition.
\end{enumerate}

\vspace{0.5in}
\subsection{The abelian category of $W(p)$-modules}
Now we denote by $W(p)$-mod the abelian category of left $W(p)$-modules. \par
Then we have the following.

\vspace{0.2in}
\noindent
{\bf Proposition 3-10.}  
The abelian category $W(p)$-mod has following properties.
\begin{enumerate}
\item The category $W(p)$-mod is Noetherian and Artinian, i.e., if $M_0 \subset M_1\subset \dots$ is an increasing sequence of objects of $W(p)$-mod then $M_n=M_{n-1}=\dots$ for some $n \geqq 0$, and if $M_0 \supseteq M_1 \supseteq \dots$ is a decreasing sequence of objects of $W(p)$-mod then $M_n=M_{n+1}=\dots$ for some $n \geqq 0$. 
\item The number of isomorphism classes of simple objects in $W(p)$-mod is finite.
\end{enumerate}

\vspace{0.2in}
\noindent
{\bf Proposition 3-11. }
For $1 \leqq s \leqq p-1$ the linear maps
$$Q_-^{[d_s^\pm]} : \mathcal{V}_s^\pm \longrightarrow \mathcal{V}_s^\mp$$
are $W(p)$-module maps.  We define $W(p)$-module $\mathcal{X}_s^\pm$ by the formulas;
$$\mathcal{X}_s^\pm=\ker Q_-^{[d_s^\mp]}(\mathcal{V}_s^\mp \longrightarrow \mathcal{V}_s^\pm).$$
Then we have the following exact sequences of $W(p)$-modules;
$$0 \longrightarrow \mathcal{X}_s^\mp \longrightarrow \mathcal{V}_s^\pm \longrightarrow \mathcal{X}_s^\pm \longrightarrow 0$$

\vspace{0.2in}
We denote $\mathcal{X}_p^\pm=\mathcal{V}_p^\pm$, $\mathcal{X}_0=\mathcal{X}_p^-$, $\mathcal{X}_p=\mathcal{X}_p^+$ where those are viewed $W(p)$-modules.  The duality on $W(p)$ is given as follows:

\vspace{0.2in}
\noindent
{\bf Proposition 3-12.}  
On $W(p)$ the following formulas hold.
\begin{enumerate}
\item $T(n)^\sigma=T(-n)$, \quad $n \in \mathbb Z$,\\
$W^a(n)^\sigma=-W^a(-n)$, \quad $a \in \{ \pm, 0 \}$, $n \in \mathbb Z$.
\item $D(\mathcal{V}_s^\pm)\simeq \mathcal{V}_s^\mp$, \quad $1\leqq s \leqq p-1$,\\
$D(\mathcal{X}_s^\pm)\simeq \mathcal{X}_s^\pm$, \quad $1\leqq s \leqq p$.
\end{enumerate}

\vspace{0.2in}
Define $\bar{X}_s^\pm \in A_0(W(p))$-mod by the following;
\begin{align*}
\bar{X}_s^+ & \equiv \mathcal{X}_s^+[h_s(0)]=C|\lambda_s(0)\rangle, \quad 1\leq s\leqq p-1, \\
\bar{X}_s^- & \equiv \mathcal{X}_s^-[h_s(1)]=C|\lambda_{-s}(0)\rangle\oplus CQ_+|\lambda_{-s}(0)\rangle, \quad 1\leqq s\leqq p-1, \\
\bar{X}_p & =\bar{X}_p^+=\mathcal{X}_p^+[h_p(0)]=C|\lambda_p(0)\rangle, \\
\bar{X}_0 & =\bar{X}_p^-=\mathcal{X}_0^-[h_0(0)]=C|\lambda_0(0)\rangle\oplus CQ_+|\lambda_0(0)\rangle .
\end{align*}

\vspace{0.2in}
\noindent
{\bf Proposition 3-13.}  
$A_0(W(p))$-modules
$$\bar{X}_s^\varepsilon : 1 \leqq s \leqq p, \varepsilon=\pm$$
are irreducible $A_0(W(p))$-modules, and all are inequivalent each other

\begin{proof}
By the definition of $h_s$, it holds that
$$h_{p-1}(1)>\dots >h_1(1)>h_0(0)>h_1(0)>\dots>h_p(0).$$
Therefore all $\bar{X}_s^\pm$ ($1\leqq s \leqq p$) are inequivalent. 
\par
For $0 \leqq s \leqq p-1$, by direct calculations we have
\begin{align*}
& W^0[0]|\lambda_{-s}(0)\rangle =\binom{-s-1}{2p-1}|\lambda_{-s}(0)\rangle, \\
& W^0[0]Q_+|\lambda_{-s}(0)\rangle =-\binom{-s-1}{2p-1}Q_-|\lambda_{-s}(0)\rangle, \\
& W^+[0]|\lambda_{-s}(0)\rangle =2\binom{-s-1}{2p-1}Q_+|\lambda_{-s}(0)\rangle, \\
& W^+[0]Q_+|\lambda_{-s}(0)\rangle =0 ,\\
& W^-[0]|\lambda_{-s}(0)\rangle =0, \\
& W^-[0]Q_+|\lambda_{-s}(0)\rangle =-\binom{-s-1}{2p-1} | \lambda_{-s}(0)\rangle.
\end{align*}

Therefore these are all irreducible $A_0(W(p))$-modules.
\end{proof}

\vspace{0.2in}
By the Proposition 3-12, we have a family of irreducible $W(p)$-modules $\mathcal{X}_s^\varepsilon$, $1 \leqq s \leqq p$ and $\varepsilon=\pm$.
\par
The structure of $\mathcal{X}_s^\pm$ as $\mathcal{L}$-modules are described in Theorem 2-7.

\vspace{0.2in}
\subsection{The structure of $W(p)$-modules $\mathcal{V}_s^\pm$, $1 \leqq s \leqq p-1$}
The structure of $W(p)$-modules $\mathcal{V}_s^\pm$ is described as follows.

\vspace{0.2in}
\noindent
{\bf Proposition 3-14.}  
For $1 \leqq s \leqq p-1$, we have:
\begin{enumerate}
\item The following equations are satisfied on $\mathcal{V}_s^+$.
\begin{enumerate}
\item $\eta_s | \lambda_{-s}(1) \rangle=Q_+ | \lambda_{-s}(0)$$=cW^+(0) | \lambda_{-s}(0) \rangle \quad (c \ne 0)$,
\item $W^-(-s) | \lambda_{-s}(1) \rangle= | \lambda_{-s}(0) \rangle$,
\item $W^0(-s) | \lambda_{-s}(1) \rangle=c' Q_+ | \lambda_{-s}(0) \rangle \quad (c' \ne 0)$.
\end{enumerate}
\item The following equations are satisfied on $\mathcal{V}_s^-$.
\begin{enumerate}
\item $\eta_s^\vee W^-(0) | \lambda_s(1) \rangle=c | \lambda_s(0) \rangle \quad (c \ne 0)$,
\item $W^-(s) | \lambda_s(1) \rangle=c'' | \lambda_s(0) \rangle \quad (c'' \ne 0)$,
\item $W^0(s)W^-(0) | \lambda_s(1) \rangle=c''' | \lambda_s(0) \rangle \quad (c''' \ne 0)$.
\end{enumerate}
\end{enumerate}

\vspace{0.5in}
\section{Construction of Log $W(p)$-modules and structure of $W(p)$-mod}
In this section, we construct $W(p)$-modules, $\mathcal{P}_s^\varepsilon$, $1\leqq s \leqq p-1$, $\varepsilon=\pm$, which we call log $W(p)$-modules, by using the logarithmic deformation of VOA $W(p)$ which is given in J. Fjeistad et al. [FFHST].  We show that the dimension of $A_0(W(p))$ is equal to $6p-1$, and give the block decomposition of $A_0(W(p))$-mod.

\vspace{0.5in}
\subsection{Construction of $\mathcal{P}_s^\pm$, $1\leqq s \leqq p-1$}
Let us fix $s$ such that $1\leqq s \leqq p-1$, and set
\begin{equation}\tag{4.1}
\mathcal{P}_s=\mathcal{V}_s^+\oplus \mathcal{V}_s^-.
\end{equation}
For each $A \in V_L$ we denote
\begin{equation}\tag{4.2}
A(z)=J^{\mathcal{V}_s^+}(A : z)\oplus J^{\mathcal{V}_s^-}(A:z) \in \mbox{End}_{\mathbb C}(\mathcal{P}_s)[[z,z^{-1}]].
\end{equation}
Then $\mathcal{P}_s$ becomes $V_L$-module by $A(z)$ for any $s$.
\par
We define operators
$$E_s^\pm(z) \in \mbox{End}_{\mathbb C}(\mathcal{P}_s)[[z,z^{-1}]]$$
by the following way;
$$E_s^\pm(z) | _{\mathcal{V}_s^\pm} = Q_-^{[d_s^\pm]}(z), \ E_s^\pm(z) | _{\mathcal{V}_s^\mp}=0.$$
Then we have
$$E_s^\pm(z) \in \mbox{Hom}_{\mathbb C}(\mathcal{V}_s^\pm, \mathcal{V}_s^\mp)[[z,z^{-1}]] \subseteq \mbox{End}_{\mathbb C}(\mathcal{P}_s)[[z,z^{-1}]].$$
\par
For each $P \in U(V)$, we denote
\begin{equation}\tag{4.3}
P = \rho^{\mathcal{V}_s^+}(P)+\rho^{\mathcal{V}_s^-}(P) \in \mbox{End}(\mathcal{P}_s).
\end{equation}
Then on End$(\mathcal{P}_s) [[ z, z^{-1} ]]$, the following properties are satisfied.  The two family of operators
\begin{equation}\tag{4.4}
\{ E_s^+(z), A(z) : A \in V_L \}, \quad \{ E_s^-(z), A(z) : A \in V_L \},
\end{equation}
are mutually local among themselves.  Also we have
\begin{equation}\tag{4.5}
E_s^+(z)E_s^+(w)=0, \quad E_s^-(z)E_s^-(w)=0.
\end{equation}
Moreover, we have
\begin{equation}\tag{4.6}
T(z)E_s^\pm(w)\sim \frac{1}{(z-w)^2}E_s^\pm(w)+ \frac{1}{(z-w)}\partial_wE_s^+(w).
\end{equation}
\par
For each $A\in V_L$, we define
\begin{align}\tag{4.7}
\Delta_s^\pm(A : z) = & (E_{s(0)}^\pm A)(z)\log z \\
& +\sum_{n \geqq 1}\frac{(-1)^n}{n}(E_{s(n)}^\pm A)(z)z^{-n} \in \mbox{End}(\mathcal{P}_s)[[z, z^{-1}]][\log z]. \notag
\end{align}
Remark that for any $A\in W(p)$, we have $(E_{s(0)}^\pm A)(z)=0$.
\par
The following two theorems can be proved easily by the methods given in [FFHST].   The construction of $W(p)$-module $\mathcal{P}_s^\pm$ is our first main result.  The analysis of the module structure of $W(p)$-modules $\mathcal{P}_s^\pm$ will be a main subjects of this paper.

\vspace{0.2in}
\noindent
{\bf Theorem 4-1.} 
There exists a unique degree preserving linear maps
\begin{equation}\tag{4.8}
\Delta_s^\pm : U(V_L) \longrightarrow \mbox{End}(\mathcal{P}_s)[\log z],
\end{equation}
$$\quad P \mapsto \Delta_s^\pm(P),$$
which satisfies the following conditions.
\begin{enumerate}
\renewcommand{\labelenumi}{(\alph{enumi})}%
\item For any $A \in V_L$, $m \in \mathbb Z$,
\begin{align}\tag{4.9}
\Delta_s^\pm(A_{(m)})=
& \left[ \int dz \ z^m(E_{s(0)}^\pm A)(z) \right] \log z \\
& +\sum_{n \geqq 1}\frac{(-1)^{n+1}}{n} \int dz \ z^{m-n}(E_{s(n)}^\pm A)(z). \notag
\end{align}
\item For all $P, Q \in U(V_L)$,
\begin{equation}\tag{4.10}
\Delta_s^\pm(PQ)=\Delta_s^\pm(P)Q+P\Delta_s^\pm(Q).
\end{equation}
\item For all $P, Q \in U$,
\begin{equation}\tag{4.11}
\Delta_s^\pm(P)\Delta_s^\pm(Q)=0.
\end{equation}
\end{enumerate}

\vspace{0.2in}
\noindent
{\bf Theorem 4-2.} 
\begin{enumerate}
\renewcommand{\labelenumi}{(\alph{enumi})}%
\item For $A \in W(p)$, define operators by
\begin{equation}\tag{4.12}
J^{\mathcal{P}_s^\pm}(A : z)=J^{\mathcal{P}_s}(A : z)+\Delta_s^\pm(A : z) \in \mbox{End}_{\mathbb C}(\mathcal{P}_s)[[z,z^{-1}]].
\end{equation}
Then these introduce a $W(p)$-module structure on $\mathcal{P}_s$ for any $s$.  We denote these $W(p)$-modules by
\begin{equation}\tag{4.13}
(\mathcal{P}_s^\pm, J^{\mathcal{P}_s^\pm}).
\end{equation}
\item We have the following exact sequence of $W(p)$-modules
$$0 \longrightarrow \mathcal{V}_s^\mp \longrightarrow \mathcal{P}_s^\pm \longrightarrow \mathcal{V}_s^\pm\longrightarrow 0.$$
\item 
\begin{equation}\tag{4.14}
J^{\mathcal{P}_s^\pm}(T, z)=T(z)+E_s^\pm(z)z,
\end{equation}
consequently we have
\begin{equation}\tag{4.15}
\rho^{\mathcal{P}_s^+}(T(0))=
\begin{cases}
T(0)+Q_-^{[p-s]}(0) & \text {on $\mathcal{V}_s^+$}, \\
T(0) & \text {on $\mathcal{V}_s^-$},
\end{cases}
\end{equation}
\begin{equation}\tag{4.16}
\rho^{\mathcal{P}_s^-}(T(0))=
\begin{cases}
T(0)+Q_-^{[s]}(0) & \text {on $\mathcal{V}_s^-$}, \\
T(0) & \text {on $\mathcal{V}_s^+$}. 
\end{cases}
\end{equation}
\end{enumerate}

\vspace{0.5in}
\subsection{Structure of $\mathcal{P}_s^\pm$, $1\leqq s \leqq p-1$}
In this subsection, we fix $s$ such that $1\leqq s \leqq p-1$.
\par
The following two theorems concerning the structures of $W(p)$-modules $\mathcal{P}_s^\pm$ are the most important results of this paper.

\vspace{0.2in}
\noindent
{\bf Theorem 4-3.}  
On the $W(p)$-module $\mathcal{P}_s^+=\mathcal{V}_s^+\oplus \mathcal{V}_s^-$, the following relations hold.
\begin{align}\tag{4.17}
& (T(0)-h_s(0)) | \lambda_{-s}(1)\rangle=| \lambda_s(0)\rangle,  \\
& (T(0)-h_s(0)) | \lambda_s(0)\rangle=0. \notag
\end{align}
\begin{align}\tag{4.18}
\eta_s | \lambda_{-s}(1)\rangle
& = \rho^{\mathcal{V}_s^+}(\eta_s) | \lambda_{-s}(1)\rangle+\Delta_s^+(\eta_s) | \lambda_{-s}(1)\rangle \\
& = Q_+ | \lambda_{-s}(0)\rangle+\Delta_s^+(\eta_s) | \lambda_{-s}(1)\rangle. \notag
\end{align}
\begin{equation}\tag{4.19}
Q_+\Delta_s^+(\eta_s) | \lambda_{-s}(1)\rangle=c | \lambda_s(1)\rangle \quad (c \ne 0),
\end{equation}
\begin{equation}\tag{4.20}
\eta_s | \lambda_s(0)\rangle= \rho^{\mathcal{V}_s^-}(\eta_s) | \lambda_s(1)\rangle=0.
\end{equation}
\begin{equation}\tag{4.21}
\eta_s^\vee | \lambda_{-s}(0)\rangle=0,
\end{equation}
\begin{align}\tag{4.22}
\eta_s^\vee W^+(0) | \lambda_{-s}(0)\rangle
& = \Delta_s^+(C_s^\vee)Q_+ | \lambda_{-s}(0)\rangle \\
& = c | \lambda_s(0)\rangle \quad (c \ne 0), \notag
\end{align}
\begin{equation}\tag{4.23}
\eta_s^\vee | \lambda_s(1)\rangle=0.
\end{equation}
\begin{align*}
\eta_s^\vee W^-(0) | \lambda_s(1)\rangle
& = \rho^{\mathcal{V}_s^-}(\eta_s^\vee) \rho^{\mathcal{V}_s^-}(W^-(0)) | \lambda_s(1)\rangle \\
& =  c | \lambda_s(0)\rangle \quad (c \ne 0).
\end{align*}

\vspace{0.2in}
\noindent
{\bf Theorem 4-4.}  
On the $W(p)$-module $\mathcal{P}_s^-=\mathcal{V}_s^+\oplus \mathcal{V}_s^-$, the following relations hold.
\begin{align}\tag{4.24}
& (T(0)-h_s(0)) | \lambda_s(0)\rangle=0, \\
& (T(0)-h_s(0)) | \lambda_{-s}(1)\rangle=0,\notag
\end{align}
\begin{align}\tag{4.25}
(T(0)-h_s(1)) | \lambda_s(1)\rangle
& = Q_+ | \lambda_{-s}(0)\rangle\\
& = cW^+(0) | \lambda_{-s}(0)\rangle \quad (c \ne 0), \notag
\end{align}
$$(T(0)-h_s(1))W^-(0) | \lambda_s(1)\rangle=c | \lambda_{-s}(0)\rangle \quad (c \ne 0).$$
\begin{align}\tag{4.26}
\eta_s | \lambda_s(0)\rangle
& = \rho^{\mathcal{V}_s^-}(\eta_s) | \lambda_s(0)\rangle+\Delta_s^-(\eta_s) | \lambda_s(0)\rangle \\
& = \Delta_s^-(\eta_s) | \lambda_s(0)\rangle \notag \\
& = c | \lambda_{-s}(0)\rangle \quad (c \ne 0), \notag
\end{align}
\begin{align}\tag{4.27}
\eta_s | \lambda_{-s}(1)\rangle
& = \rho^{\mathcal{V}_s^+}(\eta_s) | \lambda_{-s}(1)\rangle\\
& = Q_+ | \lambda_{-s}(0)\rangle. \notag
\end{align}
\begin{align*}
\eta_s^\vee | \lambda_s(1)\rangle
& = \rho^{\mathcal{V}_s^-}(\eta_s^\vee) | \lambda_s(1)\rangle+\Delta_s^-(\eta_s^\vee) | \lambda_s(1)\rangle \\
& = \Delta_s^-(\eta_s^\vee) | \lambda_s(1)\rangle \\
& = c | \lambda_{-s}(1)\rangle \quad (c \ne 0).
\end{align*}
\begin{align*}
\eta_s ^\vee | W^-(0) | \lambda_s(1)\rangle
& = \rho^{\mathcal{V}_s^-}(\eta_s^\vee)W^-(0) | \lambda_s(1)\rangle+\Delta_s^-(\eta_s^\vee)W^-(0) | \lambda_s(1)\rangle \\
& = \rho^{\mathcal{V}_s^-(\eta_s^\vee)} \rho^{\mathcal{V}_s^-}W^-(0) | \lambda_s(1)\rangle \\
& = c | \lambda_s(1)\rangle \quad (c \ne 0).
\end{align*}
$$\eta_s ^\vee | \lambda_{-s}(0)\rangle=0,$$
$$\eta_s ^\vee W^+(0) | \lambda_{-s}(0)\rangle=0.$$

\begin{proof}
We only prove Theorem 4-3.  Theorem 4-4 can be proved in the same way.
\par
In oder to prove Theorem 4-3, we express the element $\eta_s \in U(\mathcal L_{<0})[s]$ by using Bosonic operators $a(-1), a(-2), \dots$.  Consider the vector space
\begin{equation}\tag{4.28}
\mathfrak{a}_\pm=\sum_{n \geqq 1}\mathbb C a(\pm n),
\end{equation}
so that
$$\mathfrak{a}=\mathfrak{a}_+ \oplus \mathfrak{a}_-$$
These elements satisfy the following commutator relations
$$[a(m), a(n)]=m\delta_{m+n, 0}.$$
Define degree of $a(n)$ as $-n$, and consider the degreewise completed universal enveloping algebra
\begin{equation}\tag{4.29}
U(\mathfrak{a})=\sum_{d\in Z}U(\mathfrak{a})[d].
\end{equation}
Then $U(\mathfrak{a})$ is a degreewise completed tensor product of two commutative algebras such that
\begin{align}\tag{4.30}
U(\mathfrak{a}) & =U(\mathfrak{a}_-) \hat{\otimes} U(\mathfrak{a}_+), \\
U(\mathfrak{a}_\pm) & =\mathbb C[a(\pm 1), a(\pm 2), \dots ]. \notag
\end{align}
Bosonic realization of the energy-momentum tensor
$$T(z)=\frac{1}{2} : a(z)^2 : +\frac{\alpha_0}{2}\partial a(z)$$
define an algebra homomorphism
\begin{equation}\tag{4.31}
U(\mathcal{L}_{c_p}) \longrightarrow U(\mathfrak{a}) \otimes \mathbb C[a(0)].
\end{equation}
Consider a Virasoro module map
\begin{equation}\tag{4.32}
M_{h_s(0)} \longrightarrow F_{\lambda_{-s}(1)},
\end{equation}
$$\ | h_s(0\rangle \mapsto | \lambda_{-s}(1)\rangle,$$
which is $\mathbb C$-linear isomorphism up to $T(0)$ degree $h-h_s(0)\leqq s$.  Consider $\mathbb C$-linear isomorphism (4.33) and (4.34);
\begin{equation}\tag{4.33}
U(\mathcal{L}_{<0}) \longrightarrow M_{h_s(0)}
\end{equation}
$$\qquad \quad A \mapsto A | h_s(0) \rangle$$
and
\begin{equation}\tag{4.34}
U(\mathfrak{a}_-) \longrightarrow F_{\lambda_{-s}(1)}
\end{equation}
$$\qquad \quad A \mapsto A | \lambda_{-s}(1) \rangle.$$
\par
By using (4.32), (4.33) and (4.34), we get an algebra homomorphism, for $\phi_{\lambda_{-s}(1)} : U(\mathcal L_{<0}) \longrightarrow U(\mathfrak{a}_-) $
\begin{equation}\tag{4.35}
\begin{CD}
 U(\mathcal L_{<0}) @>>> U(\mathfrak{a}_-) \\
 @V{\simeq}VV
          @VV{\simeq}V \\
 M_{h_s(0)} @>>> F_{\lambda_{-s}(1)}.
\end{CD}
\end{equation}

The algebra homomorphism $\phi_{\lambda_{-s}(1)} $ is a $\mathbb C$-linear isomorphism in degree $\leqq s$.
\par
In $F_{\lambda_{-s}(1)}$, the singlar vector for $\mathcal L$ in degree $h_s(1)$ is represented by using the screening operator as
\begin{align}\tag{4.36}
Q_+ | \lambda_{-s}(0)\rangle & = \int dz \ e^{\alpha_+ \varphi(z)} | \lambda_{-s}(0)\rangle \\
& = e^{\alpha_+ \hat{a}} \int dz \ z^{s-1} e^{\alpha_+ \varphi_-(z)} | \lambda_{-s}(0)\rangle \notag \\
& = \int dz \ z^{s-1} e^{\alpha_+ \varphi_-(z)} | \lambda_{-s}(1)\rangle. \notag
\end{align}
Therefore by the map $\phi_{\lambda_{-s}(1)}$, the element $\eta_s$ is mapped to
\begin{equation}\tag{4.37}
\phi_{\lambda_{-s}(1)}(\eta_s)=\int dz \ z^{s-1} e^{\alpha_+ \varphi_-(z)} \in U(\mathfrak{a}_-) .
\end{equation}
About the element
$$\bar V_{\alpha_+}(z)=e^{\alpha_+ \varphi_-(z)}e^{\alpha_+ \varphi_+(z)} \in U(\mathfrak{a})[[z, z^{-1}]]$$
we have the following formula in $U(\mathfrak{a})$.
$$\int dz \ \bar V_{\alpha_+}(z)z^{s-1}=\int dz \ z^{s-1}e^{\alpha_+ \varphi_-(z)}+\sum_{n\geqq 1}B_na_n$$
where $B_n \in U(\mathfrak{a})[s+n]$.
\par
The map
$$\Delta_s^+ : U(\mathfrak{a}) \longrightarrow \mbox{End}(\mathcal{P}_s^+)[\log z]$$
factors through
$$\Delta_s^+ : U(\mathfrak{a}) \longrightarrow U(V_L) \longrightarrow  \mbox{End}(\mathcal{P}_s^+)[\log z]$$
and $\Delta_s^+$ is a degree preserving map which satisfies
$$\Delta_s^+(P\cdot Q)=\Delta_s^+(P)J^{P_s}(Q)+J^{P_s}(P)\Delta_s^+(Q).$$
Then we have
$$\Delta_s^+(B_n \ a_n) | \lambda_{-s}(1)\rangle=0,$$
and therefore
\begin{align*}
\Delta_s^+(\eta_s) | \lambda_{-s}(1)\rangle
& = \Delta_s^+ \left( \int dz \ z^{s-1}e^{\alpha_+ \varphi_- (z)} \right) | \lambda_{-s}(1)\rangle \\
& = \Delta_s^+ \left( \int dz \ z^{s-1} \bar V_{\alpha_+}(z)\right) | \lambda_{-s}(1)\rangle.
\end{align*}
Now we see
\begin{align*}
& \langle \lambda_s(1) | Q_+\Delta_s^+(\eta_s) | \lambda_{-s}(1)\rangle \\
& = \int dz \ z^{s-1}\langle \lambda_s(1) | Q_+\Delta_s^+(\bar V_{\alpha_+}(z)) | \lambda_{-s}(1)\rangle \\
& = \int dy \int dw \int dz \ z^{s-1}\langle \lambda_s(1) | V_{\alpha_+}(y)Q_+^{[p-s]}(w)\bar V_{\alpha_+}(z) | \lambda_{-s}(1)\rangle \\
& = \int dy \int dw_1 \cdots dw_{p-s} \int  dz \ z^{s-1}\langle \lambda_s(1) | V_{\alpha_+}(y) V_{\alpha_-}(w_1) \\
& \qquad \qquad \qquad \qquad \qquad \qquad \qquad \dots V_{\alpha_-}(w_{p-s}) \bar V_{\alpha_+}(z)  | \lambda_{-s}(1)\rangle \\
& \ne 0
\end{align*}
Consequently we get
$$Q_+\Delta_s^+(\eta_s) | \lambda_{-s}(1)\rangle=\mbox{const.} | \lambda_s(1)\rangle \quad (\mbox{const.} \ne 0).$$
\par
Proof of $\Delta_s^+(C_s^V)Q_+ | \lambda_{-s}(0)\rangle=\mbox{const.} | \lambda_s(0)\rangle$ (const.$\ \ne 0$) can be done exactly in the same way.
\end{proof}

By using the results of \S 4-2, it is easy to verify the following structure of projective module $\mathcal{P}_s^\pm$.

\vspace{0.2in}
\noindent
{\bf Proposition 4-5.} 
\begin{enumerate}
\item The socles sequence of $\mathcal{P}_s^+$ is
$$S_1(\mathcal{P}_s^+) \subseteq S_2(\mathcal{P}_s^+) \subseteq S_3(\mathcal{P}_s^+)=\mathcal{P}_s^+,$$
$$S_1(\mathcal{P}_s^+) \simeq \mathcal{X}_s^+ , \ S_2(\mathcal{P}_s^+) / S_1(\mathcal{P}_s^+) \simeq \mathcal{X}_s^- \oplus \mathcal{X}_s^- , \ S_3(\mathcal{P}_s^+) / S_2(\mathcal{P}_s^+) \simeq \mathcal{X}_s^+.$$
\item The socles sequence of $\mathcal{P}_s^-$ has following structures
$$S_1(\mathcal{P}_s^-) \subseteq S_2(\mathcal{P}_s^-) \subseteq S_3(\mathcal{P}_s^-)=\mathcal{P}_s^-,$$
$$S_1(\mathcal{P}_s^-) \simeq \mathcal{X}_s^- , \ S_2(\mathcal{P}_s^-) / S_1(\mathcal{P}_s^-) \simeq \mathcal{X}_s^+ \oplus \mathcal{X}_s^+ , \ S_3(\mathcal{P}_s^-) / S_2(\mathcal{P}_s^-) \simeq \mathcal{X}_s^-.$$
\end{enumerate}

\vspace{0.5in}
\subsection{Structure of $A_0(W(p))$} 
Let us consider the $W(p)$-module $\mathcal{P}_s^+$, $1 \leqq s \leqq p-1$.  Define
$$\bar{P}_s^+=\mathcal{P}_s^+[h_s(0)]=\mathbb C| \lambda_s(0)\rangle + \mathbb C| \lambda_{-s}(0)\rangle$$
Then $\bar{P}_s^+$ is a $A_0(W(p))$-module.  By Theorem 4-3 we see
$$\rho_s^{\mathcal{P}_s^+}(T(0)-h_s(0)) | \lambda_{-s}(1)\rangle=| \lambda_s(0)\rangle ,$$
$$\rho_s^{\mathcal{P}_s^+}(T(0)-h_s(0)) | \lambda_s(0)\rangle=0 .$$
\par
For each $1 \leqq s \leqq p$ and $\varepsilon=\pm$, we define finite dimensional algebra $I_s^\varepsilon$ by the following way.
\begin{enumerate}
\item Case 1.  $1 \leqq s \leqq p-1$, $\varepsilon=1$. \\
Consider the algebra homomorphism
$$\rho_s^+ : A_0(W(p)) \longrightarrow \mbox{End}(\bar{P}_s^+)=M_2(\mathbb C).$$
Then $\mbox{Image} \ (\rho_s^+) \subseteq M_2(\mathbb C)$ contain two dimensional algebra
$$I_s^+=\left\{  \left(
\begin{array}{cc}
a& 0\\
b& a
\end{array}
\right) : a, b \in \mathbb C
\right\}$$

\item Case 2.  $s=p$, $\varepsilon=\pm$ or $1 \leqq s \leqq p-1$, $\varepsilon=-1$. \\
Consider the algebra homomorphism for any $s$
$$\rho_s^\varepsilon : A(W(p)) \longrightarrow \mbox{End} (\bar{X}_s^\varepsilon).$$
Since $\bar{X}_s^\varepsilon$ is an irreducible $A_0(W(p))$-module, the map $\rho_s^\varepsilon$ is surjective.  Set
$$I_s^\varepsilon=\mbox{Image} ( \rho_s^\varepsilon) =\mbox{End} (\bar{X}_s^\varepsilon).$$
\end{enumerate}

\vspace{0.2in}
\noindent
{\bf Theorem 4-6.} 
\begin{enumerate}
\item The algebra $A_0(W(p))$ is isomorphic to $I=\displaystyle\sum_{s=1}^p \sum_{\varepsilon=\pm} I_s^\varepsilon.$
\item The center of $A_0(W(p))$ is generated by $[T(0)] .$
\item The set of inequivalent irreducible $A_0(W(p))$-modules is \{ $\bar{X}_s^\varepsilon$, $1\leqq s \leqq p$, $\varepsilon=\pm$\}.
\end{enumerate}

\begin{proof}
We see $\dim_{\mathbb C} I=6p-1$ and $\dim_{\mathbb C}A_0(p) \leqq 6p-1$.  On the other hand, by definition we see $\dim_{\mathbb C} I \leqq \dim_{\mathbb C} A_0(p)$.  So we get the proposition.
\end{proof}

\vspace{0.5in}
\subsection{$A_0(W(p))$-mod}
For each $s=1, \dots, p$ and $\varepsilon=\pm$, we define the full abelian subcategories $\bar{\mathcal C}_s^\varepsilon$ of $A_0(W(p))$-mod such that an element $M\in A_0(W(p))$-mod belongs to $\bar{\mathcal C}_s^\varepsilon$, if and only if the irreducible components of $M$ are $\bar{X}_s^\varepsilon$.  Then we have the following theorem from Theorem 4-6.

\vspace{0.2in}
\noindent
{\bf Theorem 4-7.}
\begin{enumerate}
\item The abelian category $A_0(W(p))$-mod has the following block decomposition
$$A_0(W(p))\mbox{-mod}=\sum_{s=1}^p \sum_{\varepsilon=\pm} \bar{\mathcal C}_s^\varepsilon .$$
\item For $s=p$, $\varepsilon=\pm$ or $1\leqq s \leqq p-1$, $\varepsilon=-$, the abelian category $\bar{\mathcal C}_s^\varepsilon$ is semi-simple with a simple object $\bar{X}_s^\varepsilon$ .
\item For $1\leqq s \leqq p-1$, the set of indecomposable objects in the abelian category $\bar{\mathcal C}_s^+$ is \{$\bar{X}_s^+$, $\bar{P}_s^+$\}.  Moreover, we have the non-trivial exact sequence of $A_0(W(p))$-mod
$$0 \longrightarrow \bar{X}_s^+ \longrightarrow \bar{P}_s^+ \longrightarrow \bar{X}_s^+ \longrightarrow 0.$$
\end{enumerate}

\vspace{0.5in}
\subsection{Block decomposition of the abelian category $W(p)$-mod}
The following Theorem 4-8 is proved in [AM1], [AM2]

\vspace{0.2in}
\noindent
{\bf Theorem 4-8.}  
$\mbox{Ext}_{W(p)}^1(\mathcal{X}_{s_1}^{\varepsilon_1}, \mathcal{X}_{s_2}^{\varepsilon_2})=0$
for $1 \leqq s_1 \ne s_2 \leqq p$ and $\varepsilon_1, \varepsilon_2=\pm$.

\vspace{0.2in}
For each $0 \leqq s \leqq p$ we denote by $C_s$ the full abelian category of $W(p)$-mod such that $M \in W(p)$-mod belong to $C_s$ if and only if $M$ has Jordan-H\"{o}lder sequence whose factors are $\mathcal{X}_s^\pm$ if $1\leqq s \leqq p-1$, and $\mathcal{X}_s$ if $s=0$ or $p$, respectively.
\par
Then by virtue of theorem 4-8, we have the following.

\vspace{0.2in}
\noindent
{\bf Theorem 4-9.}  
The abelian category $W(p)$-mod has the following block decomposition
$$W(p)\mbox{-mod}=\sum_{s=0}^p C_s ,$$
with the properties:
\begin{enumerate}
\item Each element of  $W(p)$-mod has the unique decomposition
$$M=\sum_{s=0}^p C_s \quad \mbox{with} \quad M_s \in C_s.$$
\item For any $M \in C_s$, $N \in C_{s'}$,
$$\mbox{Ext}^{\bullet}(M, N)=0 \quad \mbox{if} \quad s \ne s'.$$
\end{enumerate}

\vspace{0.2in}
\noindent
{\bf Proposition 4-10. }
For each $0\leqq s \leqq p$, any element $M \in C_s$ has following eigenspace decompositions
$$M=\sum_{h\in H_s} M[h],$$
where $M[h]=\{ m \in M : (T(0)-h)^n m=0$ for some $n \geqq 1$\}, and $\dim_{\mathbb C} M[h]<\infty$ for all $h \in \mathbb C$.
\par
\vspace{0.2in}
The following Proposition is very important in this paper.

\vspace{0.2in}
\noindent
{\bf Proposition 4-11.}  
Let $s$ be an integer such that $1 \leqq s \leqq p-1$, and let $M, N \in C_s$, and $f : M \rightarrow N$ be a $W(p)$-module map.  If $f$ is a vector space isomorphism of degree $h$, for $h-h_s(0) \leqq s$, then $f$ is a $W(p)$-module isomorphism.

\begin{proof}
Category $C_s$ has simple objects $\mathcal{X}_s^\pm$, and the highest weight of $\mathcal{X}_s^+$ and $\mathcal{X}_s^-$ are $h_s(0)$ and $h_s(1)$, respectively.  Note that $h_s(1)-h_s(0)=s$.
\par
Consider the kernel and the cokernel of $f$, then by the condition of $\mathcal{P}$ the weight $h$ satisfies $h-h_s(0)>s$.  This shows that the kernel and the cokernel of $f$ must be zero.
\end{proof}

\vspace{0.2in}
\noindent
{\bf Proposition 4-12. } 
\begin{enumerate}
\item For $s=p$ and $\varepsilon=\pm$ or $1\leqq s \leqq p-1$, $\varepsilon=-1$, we have isomorphism of $W(p)$-modules
$$U(W(p))\underset{F_0(W(p))}{\otimes}\bar{X}_s^\varepsilon \simeq \mathcal{X}_s^\varepsilon$$
\item For $1\leqq s \leqq p-1$, $\varepsilon=\pm$, an element $M$ of $C_s$ is a direct sum of $\mathcal{X}_s^\varepsilon$ if and only if $M=\displaystyle\sum_{h\in H_s^\varepsilon}M[h]$
\end{enumerate}

\vspace{0.5in}
\section{Projectivity of $\mathcal{P}_s^\pm$}
In this section we show that $\mathcal{P}_s^\pm$, $1\leqq s \leqq p$, are projective covers of simple modules $\mathcal{X}_s^\pm$, $1\leqq s \leqq p$.

\vspace{0.5in}
\subsection{The structure of $\mbox{Ext}^1(\mathcal{X}_s^\varepsilon, \mathcal{X}_{s'}^{\varepsilon'})$}
The following two theorems are part of our main results.

\vspace{0.2in}
\noindent
{\bf Theorem 5-1.}  
For $1\leqq s \leqq p$, $\varepsilon=\pm$
\begin{equation}\tag{5.1}
\mbox{Ext}^1(\mathcal{X}_s^\varepsilon, \mathcal{X}_s^\varepsilon)=0
\end{equation}

\begin{proof}
We divide a proof into two cases.
\begin{enumerate}
\item Case 1. $s=p$ and $\varepsilon=\pm$, or $1\leqq s\leqq p-1$, $\varepsilon=-1$. \\
We denote $X=\mathcal{X}_s^\varepsilon$ for simplicity, and consider exact sequence of $W(p)$-modules
$$0 \longrightarrow X \longrightarrow E \longrightarrow X \longrightarrow 0.$$
The lowest $T(0)$ degree part of this exact sequence gives an exact sequence of $A_0(W(p))$-modules
$$0\longrightarrow \bar{X} \longrightarrow \bar{E}\longrightarrow \bar{X} \longrightarrow 0.$$
Then this exact sequence belongs to the block $\bar{C}_s^\varepsilon$.  This case $\bar{C}_s^\varepsilon$ is a semi-simple category whose simple object is $X$.  So we get $E = X \oplus X$ as $W(p)$-module by proposition 4-11(1).\\

\item Case 2.  $1\leqq s\leqq p-1$, $\varepsilon=+$\\
We denote $X=\mathcal{X}_s^+$ and consider an exact sequence of $W(p)$-modules
$$0\longrightarrow X\longrightarrow E \longrightarrow X\longrightarrow 0.$$
Then the lowest $T(0)$ degree part of this sequence gives an exact sequence of $A_0(W(p))$-modules
$$0 \longrightarrow \bar{X} \longrightarrow \bar{E} \longrightarrow \bar{X }\longrightarrow 0.$$
This sequence belongs to the block $\bar{C}_s^+$, $1\leqq s\leqq p-1$, and $E=X \oplus X$ as $W(p)$-modules if and only if $\bar{E}=\bar{X} \oplus \bar{X}$ in $\bar{C}_s^+$, that is, $[T(0)]$ acts on $\bar{E}$ semi-simple.  Consider $\mathcal{L}$-module $M=U(\mathcal{L})(\bar{E})\subseteq E$.  Then as $\mathcal{L}$-modules it has following exact sequences
$$0\longrightarrow L_{h_s(0)} \longrightarrow M \longrightarrow L_{h_s(0)} \longrightarrow 0$$
Then by proposition 2-11 for Virasoro modules gives $\bar{E}=\bar{X} \oplus \bar{X}$ as $A_0(W(p))$-modules.  Thus by Proposition 4-11(1) we have 
$$E=X \oplus X$$
as $W(p)$-modules
\end{enumerate}
\end{proof}

\vspace{0.2in}
We define the $W(p)$-module $\mathcal{Y}_s^+$, $1 \leqq s \leqq p-1$, by the following exact sequence;
\begin{equation}\tag{5.2}
0 \longrightarrow \mathcal{X}_s^+ \longrightarrow \mathcal{P}_s^+ \longrightarrow \mathcal{Y}_s^+\longrightarrow 0.
\end{equation}

\vspace{0.2in}
\noindent
{\bf Theorem 5-2.}  
For $1\leqq s \leqq p-1$ we have
$$\mbox{Ext}^1(\mathcal{X}_s^\pm, \mathcal{X}_s^\mp)= \mathbb C^2.$$

\begin{proof}
We first prepare some notations.  By the duality in $W(p)$-mod, we have $D(\mathcal{X}_s^\varepsilon)\cong \mathcal{X}_s^\varepsilon$, thus it is suficient to prove $\mbox{Ext}^1(\mathcal{X}_s^+, \mathcal{X}_s^-)\simeq \mathbb C^2$.
\par
First we fix the following elements of $\mathcal{X}_s^\pm$;
$$u=|\lambda_s(0)\rangle\in \mathcal{X}_s^+ \subseteq \mathcal{V}_s^- ,$$
$$v_-=|\lambda_{-s}(0)\rangle,\quad v_+=Q_+|\lambda_{-s}(0)\rangle\in \mathcal{X}_s^- \subseteq \mathcal{V}_s^+ .$$
The element
$$v_+=Q_+ | \lambda_{-s}(0)\rangle \in F_{\lambda_{-s}(1)}[h_s(1)]$$
is a singlar vector of Virasoro module $F_{\lambda-s(1)}[h_s(1)]$, and the sequence of $\mathcal{L}_{c_p}$-module maps;
$$M_{h_s(2)} \longrightarrow M_{h_s(0)} \longrightarrow F_{\lambda_{-s}(1)}$$
is exact, so we have
$$v_+=\eta_s | \lambda_{-s}(1)\rangle.$$
Now we give a proof of Theorem 5-2.
\par
Consider an exact sequence and the elements in $E$, and fix element $\tilde{u} \in E$,
$$0 \longrightarrow \mathcal{X}_s^- \longrightarrow E \longrightarrow \mathcal{X}_s^+\longrightarrow 0$$
$$\mathcal{X}_s^- \ni v_+, v_-, \qquad E \ni \tilde u \mapsto u$$
such that elements $\tilde u \in E[h_s(0)]$ is mapped to $u$ in $\mathcal{X}_s^+$, which is uniquely determined.  Set
$$\eta_s\tilde u=a_+v_++a_-v_- ,$$
Note that if $a_+=0$ and $a_-=0$, then $[E]=0$ in $\mbox{Ext}^1(\mathcal{X}_s^+, \mathcal{X}_s^-)$.
\par
By the definition of $W(p)$-module $\mathcal{P}_s^+$ we see the exact sequence
$$0 \longrightarrow \mathcal{V}_s^- \longrightarrow \mathcal{P}_s^+ \longrightarrow \mathcal{V}_s^+\longrightarrow 0.$$
We define two $W(p)$-submodules $E_1$ and $E_2$ of $\mathcal{Y}_s$ by
\begin{align*}
& E_1=\mathcal{V}_s^- / \mathcal{X}_s^+ \hookrightarrow \mathcal{Y}_s^+=\mathcal{P}_s^+ / \mathcal{X}_s^+ , \\
& E_2=U(W(p))(\mathbb C |\lambda_{-s}(0)\rangle \oplus \mathbb CQ_+ |\lambda_{-s}(0)\rangle) \subseteq \mathcal{Y}_s^+.
\end{align*}
Then by using Theorem 4-4, it is easy to show that the $W(p)$-modules $E_1$ and $E_2$ are both isomorphic to the $W(p)$-module $\mathcal{X}_s^+$.
\par
Consequently $W(p)$-module $\mathcal{Y}_s / E_1$ is canonically isomorphic to $\mathcal{V}_s^+$.  Let us introduce a $W(p)$-module $\mathcal{Y}_s^+ / E_2=\mathcal{V}_s^{+\vee}$.  Then we have exact sequences
$$0\longrightarrow \mathcal{X}_s^-\longrightarrow \mathcal{V}_s^+\longrightarrow \mathcal{X}_s^+\longrightarrow 0,$$
$$0\longrightarrow \mathcal{X}_s^-\longrightarrow \mathcal{V}_s^{+\vee}\longrightarrow \mathcal{X}_s^+\longrightarrow 0.$$
For $\mathcal{V}_s^+$, we have $a_+=1$ and $a_-=0$.  Also for $\mathcal{V}_s^{+\vee}$, we have $a_+=0$ and $a_-=1$.  Consequently $[\mathcal{V}_s^+]$ and $[\mathcal{V}_s^-]$ are linearly independent in $\mbox{Ext}^1(\mathcal{X}_s^+, \mathcal{X}_s^-)$.
\end{proof}

\vspace{0.2in}
\noindent
{\bf Proposition 5-3. }
The subcategories $C_0$ and $C_p$ of $W(p)$-mod are semi-simple with only one simple object $\mathcal{X}_0$, and $\mathcal{X}_p$, respectively.

\begin{proof}
By the theorem 5-2 we have $\mbox{Ext}^1(\mathcal{X}_0, \mathcal{X}_0)=0$, $\mbox{Ext}^1(\mathcal{X}_p, \mathcal{X}_p)=0$.  Therefore we have proved the statement.
\end{proof}

We denote $\mathcal{P}_p=\mathcal{P}_p^+ = \mathcal{V}_p=\mathcal{X}_p$ and $\mathcal{P}_0=\mathcal{P}_p^-=\mathcal{V}_0=\mathcal{X}_0$.  Then these two modules are projective modules in $W(p)$-mod.

\vspace{0.2in}
\subsection{Projectivity of $\mathcal{P}_s^+$, $1\leqq s \leqq p-1$}
We fix $s$ such that $1\leqq s \leqq p-1$.

\vspace{0.2in}
\noindent
{\bf Proposition 5-4.} 
One has
\begin{equation}\tag{5.3}
\mbox{Hom}(\mathcal{V}_s^+, \mathcal{X}_s^+)=\mathbb C, \quad \mbox{Hom}(\mathcal{V}_s^{+\vee}, \mathcal{X}_s^+)\simeq \mathbb C,
\end{equation}
$$\mbox{Hom}(\mathcal{V}_s^+, \mathcal{X}_s^-)=0, \quad \mbox{Hom}(\mathcal{V}_s^{+\vee}, \mathcal{X}_s^-)=0,$$
$$\mbox{Hom}(\mathcal{Y}_s^+, \mathcal{X}_s^+)=\mathbb C, \quad \mbox{Hom}(\mathcal{Y}_s^+, \mathcal{X}_s^-)=0,$$
\begin{equation}\tag{5.4}
\mbox{Ext}^1(\mathcal{V}_s^+, \mathcal{X}_s^-)\simeq \mathbb C,
\end{equation}
$$\mbox{Ext}^1(\mathcal{V}_s^{+\vee}, \mathcal{X}_s^-)\simeq \mathbb C,$$
\begin{equation}\tag{5.5}
\mbox{Ext}^1(\mathcal{Y}_s^+, \mathcal{X}_s^-)=0.
\end{equation}

\begin{proof}
These follow from the definitions and results obtained in \S 5-1.  So we only prove (5.5).  Consider exact sequence
$$0 \longrightarrow (\mathcal{X}_s^-)^2 \longrightarrow \mathcal{Y}_s^+\longrightarrow \mathcal{X}_s^+\longrightarrow 0,$$
which give an exact sequence 
$$0 \longrightarrow \mbox{Hom}((\mathcal{X}_s^-)^2, \mathcal{X}_s^-) \longrightarrow \mbox{Ext}^1(\mathcal{X}_s^+, \mathcal{X}_s^-) \longrightarrow \mbox{Ext}^1(\mathcal{Y}_s^+, \mathcal{X}_s^-)\longrightarrow 0.$$
Therefore as discussed in \S 5-1, we have
$$\mbox{Ext}^1(\mathcal{Y}_s^+, \mathcal{X}_s^-)=0 .$$
\end{proof}

\vspace{0.2in}
\noindent
{\bf Proposition 5-5. } 
The linear map
\begin{equation}\tag{5.6}
U(W(p)) \underset{F_0(U(W(p)))}{\otimes}\bar{X}_s^+ \longrightarrow \mathcal{Y}_s^+
\end{equation}
is an isomorphism of $W(p)$-modules.

\begin{proof}
Consider canonical map
$$U(W(p)) \underset{F_0(U(W(p)))}{\otimes}\bar{X}_s^+\longrightarrow \mathcal{Y}_s^+$$
This map is surjective, and the kernel is isomorphic to $(\mathcal{X}_s^-)^l$ for some $l \geqq 0$.  But $\mbox{Ext}^1(\mathcal{Y}_s^+, \mathcal{X}_s^-)=0$.  So we have $l=0$
\end{proof}

\vspace{0.2in}
\noindent
{\bf Proposition 5-6. } 
\begin{equation}\tag{5.7}
\mbox{Ext}^1(\mathcal{Y}_s^+, \mathcal{X}_s^+)\simeq \mathbb C
\end{equation}

\begin{proof}
Since $[\mathcal{P}_s^+]\ne 0$ in $\mbox{Ext}^1(\mathcal{Y}_s^+, \mathcal{X}_s^+)$, we have $\dim_{\mathbb C}\mbox{Ext}^1(\mathcal{Y}_s^+, \mathcal{X}_s^+) \geqq 1$.  Consider an exact sequence of $W(p)$-modules
$$0 \longrightarrow \mathcal{X}_s^+\longrightarrow E \longrightarrow \mathcal{Y}_s^+ \longrightarrow 0,$$
and fix elements $u_0 \in \mathcal{X}_s^+[h_s(0)]\simeq \mathbb C$ and $u_1 \in \mathcal{X}_s^+[h_s(0)]\simeq \mathbb C$. Take an element $\tilde{u}_0\in E[h_s(0)]\simeq {\mathbb C}^2$ which is mapped to $u_0$ in $\mathcal{Y}_s^+$.  Then we have $(T(0)-h_s(0))\tilde{u}=cu$ for some $c\in \mathbb C$.   We assume $c=0$, then by Proposition 5-4 we have a following $W(p)$-module map
$$\mathcal{Y}_s^+=U(W(p)) \underset{F_0(U(W(p)))}{\otimes}\bar{X}_s^+\longrightarrow E,$$
which is a lifting of $E \rightarrow \mathcal{Y}_s^+ \rightarrow 0$.  Thus $[E]=0$ in $\mbox{Ext}^1(\mathcal{Y}_s^+, \mathcal{X}_s^+)$.  This show that $\dim_{\mathbb C}\mbox{Ext}^1(\mathcal{Y}_s^+, \mathcal{X}_s^+)\leqq 1$.   Therefore we get the result.
\end{proof}

\vspace{0.2in}
\noindent
{\bf Proposition 5-7.} 
\begin{equation}\tag{5.8}
\mbox{Ext}^1(\mathcal{P}_s^+, \mathcal{X}_s^+)=0
\end{equation}

\begin{proof}
This follows easily by the following exact sequence and Proposition 5-6
$$0 \longrightarrow \mathcal{X}_s^+ \longrightarrow \mathcal{P}_s^+ \longrightarrow \mathcal{Y}_s^+ \longrightarrow 0.$$
\end{proof}

\vspace{0.2in}
\noindent
{\bf Proposition 5-8. } 
For any $s$ we have
\begin{equation}\tag{5.9}
\mbox{Ext}^1(\mathcal{P}_s^+, \mathcal{X}_s^-)=0.
\end{equation}

\begin{proof}
Consider an exact sequence
$$0 \longrightarrow \mathcal{X}_s^+\longrightarrow \mathcal{P}_s^+\longrightarrow \mathcal{Y}_s^+\longrightarrow 0 .$$
Since $\mbox{Ext}^1(\mathcal{Y}_s^+, \mathcal{X}_s^-)=0$ this gives
\begin{equation}\tag{5.10}
0 \longrightarrow \mbox{Ext}^1(\mathcal{P}_s^+, \mathcal{X}_s^-) \longrightarrow \mbox{Ext}^1(\mathcal{X}_s^+, \mathcal{X}_s^-).
\end{equation}
Let us consider an exact sequence
\begin{equation}\tag{5.11}
0 \longrightarrow \mathcal{X}_s^- \longrightarrow E \longrightarrow \mathcal{P}_s^+ \longrightarrow 0,
\end{equation}
and define $\bar u_1= | \lambda_s(0) \rangle$, $\bar u_0=| \lambda_{-s}(0)\rangle$ in $\mathcal{P}_s^+$.  Take the elements $u_i \in E$ which are mapped to $\bar u_i \in \mathcal{P}_s^+$ for $i=0,1$.
\par
Then we have $(T(0)-h_1(0))u_0=u_1$, $(T(0)-h_1(0))u_1=0$ and $T(n)u_i=0$, for $n \geqq 1$, $i=0,1$.  By Proposition 2-11, we have $\eta_s(u_1)=0$.  This shows that $[E]=0$ in $\mbox{Ext}^1(\mathcal{X}_s^+, \mathcal{X}_s^-)$.  Consequently $[E]=0$ in $\mbox{Ext}^1(\mathcal{P}_s^+, \mathcal{X}_s^-)$.
\end{proof}

\vspace{0.2in}
\noindent
{\bf Theorem 5-9.}  
\begin{enumerate}
\item $\mathcal{P}_a^+$ are projective $W(p)$-modules.
\item For all $s$, $\mathcal{P}_s^+\rightarrow \mathcal{X}_s^+\rightarrow 0$ are projective covers.
\end{enumerate}

\begin{proof}
These are direct consequences of Theorem 5-7 and Theorem 5-8.
\end{proof}

\vspace{0.2in}
\noindent
{\bf Proposition 5-10.} 
\begin{enumerate}
\item $\mbox{Ext}^1(\mathcal{V}_s^+, \mathcal{V}_s^-)\simeq \mathbb C$, $\mbox{Ext}^1(\mathcal{V}_s^{+\vee}, \mathcal{V}_s^{-\vee})\simeq \mathbb C$.
\item These two vector spaces in (1) have generators $[\mathcal{P}_s^+]$.
\end{enumerate}

\begin{proof}
Consider exact sequences
$$0 \longrightarrow \mathcal{V}_s^-\longrightarrow \mathcal{P}_s^+ \longrightarrow \mathcal{V}_s^+\longrightarrow 0,$$
$$0 \longrightarrow \mathcal{V}_s^{-\vee}\longrightarrow \mathcal{P}_s^+ \longrightarrow \mathcal{V}_s^{+\vee}\longrightarrow 0.$$
\par
Then statements follow immediately.
\end{proof}

\vspace{0.2in}
\noindent
{\bf Proposition 5-11.} 
One has
\begin{equation}\tag{5.12}
\mbox{Ext}^1(\mathcal{V}_s^+, \mathcal{V}_s^{-\vee})=0,
\end{equation}
$$\mbox{Ext}^1(\mathcal{V}_s^{+\vee}, \mathcal{V}_s^-)=0.$$

\begin{proof}
A proof is same as the one of Proposition 5-10. 
\end{proof}

\vspace{0.2in}
\noindent
{\bf Proposition 5-12. } 
\begin{enumerate}
\item $D(\mathcal{P}_s^+)\simeq \mathcal{P}_s^+$.
\item $\mathcal{P}_s^+$ is an injective module.
\end{enumerate}

\begin{proof}
(2) follows from (1), since $\mathcal{P}_s^+$ is a generator of $\mbox{Ext}^1(\mathcal{V}_s^+, \mathcal{V}_s^-)\simeq \mathbb C$.  But $D(\mathcal{V}_s^\pm) \simeq \mathcal{V}_s^\mp$, and then we have $D(\mathcal{P}_s^+)\simeq \mathcal{P}_s^+$.
\end{proof}

\vspace{0.2in}
\noindent
{\bf Proposition 5-13.} 
\begin{equation}\tag{5.13}
\mbox{Ext}^1(\mathcal{V}_s^+, \mathcal{X}_s^+)\simeq 0,
\end{equation}
$$\mbox{Ext}^1(\mathcal{X}_s^{+\vee}, \mathcal{X}_s^+)\simeq 0.$$

\begin{proof}
Consider the exact sequence
$$0\longrightarrow \mathcal{V}_s^-\longrightarrow \mathcal{P}_s^+\longrightarrow \mathcal{V}_s^+\longrightarrow 0.$$
Then we have the exact sequence
$$0=\mbox{Hom}(\mathcal{V}_s^-, \mathcal{X}_s^+)\overset{\sim}{\longrightarrow} \mbox{Ext}^1(\mathcal{V}_s^+, \mathcal{X}_s^+).$$
Then we can prove $\mbox{Ext}^1(\mathcal{V}_s^+, \mathcal{X}_s^+)\simeq 0$ similary.
\end{proof}

\vspace{0.5in}
\subsection{Projectivity of $\mathcal{P}_s^-$, $1\leqq s\leqq p-1$}
We fix $1\leqq s \leqq p-1$, and define the $W(p)$-module $\mathcal{Y}_s^-$ by the exact sequence
\begin{equation}\tag{5.14}
0 \longrightarrow \mathcal{X}_s^- \longrightarrow \mathcal{P}_s^- \longrightarrow \mathcal{Y}_s^-\longrightarrow 0.
\end{equation}

\vspace{0.2in}
\noindent
{\bf Proposition 5-14.}  
We have
\begin{enumerate}
\item $\mbox{Ext}^1(\mathcal{V}_s^-, \mathcal{X}_s^+)\simeq \mathbb C$, \quad $\mbox{Ext}^1(\mathcal{V}_s^{-\vee}, \mathcal{X}_s^+)\simeq \mathbb C$,
\item $\mbox{Ext}^1(\mathcal{Y}_s^-, \mathcal{X}_s^+)\simeq 0 .$
\end{enumerate}

\begin{proof}
Consider an exact sequence
$$0\longrightarrow \mathcal{X}_s^+\longrightarrow \mathcal{V}_s^-\longrightarrow \mathcal{X}_s^-\longrightarrow 0.$$
This gives exact sequence
$$0\longrightarrow \mbox{Hom}(\mathcal{X}_s^+, \mathcal{X}_s^+) \longrightarrow \mbox{Ext}^1(\mathcal{X}_s^-, \mathcal{X}_s^+)\longrightarrow \mbox{Ext}^1(\mathcal{V}_s^-, \mathcal{X}_s^+)\longrightarrow 0,$$
and then we get $\mbox{Ext}^1(\mathcal{V}_s^-, \mathcal{X}_s^+)\simeq \mathbb C$.
\par
In the same way, we can conclude $\mbox{Ext}^1(\mathcal{V}_s^{-\vee}, \mathcal{X}_s^+)\simeq \mathbb C$
\par
Consider the exact sequence
$$0 \longrightarrow \mathcal{X}_s^+\longrightarrow \mathcal{Y}_s^-\longrightarrow \mathcal{V}_s^-\longrightarrow 0.$$
Then we have an exact sequence
$$0 \longrightarrow \mbox{Hom}(\mathcal{X}_s^+, \mathcal{X}_s^+) \longrightarrow \mbox{Ext}^1(\mathcal{V}_s^-, \mathcal{X}_s^+)\longrightarrow \mbox{Ext}^1(\mathcal{Y}_s^-, \mathcal{X}_s^+)\longrightarrow 0.$$
The statement (2) follows from this sequence.
\end{proof}

\vspace{0.2in}
\noindent
{\bf Proposition 5-15.}  
$$\mbox{Ext}^1(\mathcal{P}_s^-, \mathcal{X}_s^+)=0.$$

\begin{proof}
Consider the exact sequence
$$0 \longrightarrow \mathcal{V}_s^+\longrightarrow \mathcal{P}_s^-\longrightarrow \mathcal{V}_s^-\longrightarrow 0.$$
Then we have exact sequence
$$0 \longrightarrow \mbox{Hom}(\mathcal{V}_s^+, \mathcal{X}_s^+) \longrightarrow \mbox{Ext}^1(\mathcal{V}_s^-, \mathcal{X}_s^+)\longrightarrow \mbox{Ext}^1(\mathcal{P}_s^-, \mathcal{X}_s^+)\longrightarrow \mbox{Ext}^1(\mathcal{V}_s^+, \mathcal{X}_s^+).$$
 By Proposition 5-14, we have $\mbox{Ext}^1(\mathcal{V}_s^-, \mathcal{X}_s^+)\simeq \mathbb C$.  Therefore by Proposition 5-13 we have $\mbox{Ext}^1(\mathcal{V}_s^+, \mathcal{X}_s^+)=0$.
\end{proof}

\vspace{0.2in}
\noindent
{\bf Proposition 5-16.} 
\begin{enumerate}
\item $\mbox{Ext}^1(\mathcal{Y}_s^-, \mathcal{X}_s^-) \simeq \mathbb C.$
\item The element $[\mathcal{P}_s^-]$ is a generator of $\mbox{Ext}^1(\mathcal{Y}_s^-, \mathcal{X}_s^-)$.
\end{enumerate}

\begin{proof}
Since the element $[\mathcal{P}_s^-]$ is non-zero element of $\mbox{Ext}^1(\mathcal{Y}_s^-, \mathcal{X}_s^-)$, it is surficient to prove $\dim_{\mathbb C}\mbox{Ext}^1(\mathcal{Y}_s^-, \mathcal{X}_s^-) \leqq 1$
\par
We fix elements of $\mathcal{Y}_s^-=\mathcal{P}_s^- / \mathcal{X}_s^-$ by the following way.
\begin{align*}
& v_+= | \lambda_s(1) \rangle, \ v_-=W^-(0)v_+ \in \mathcal{Y}_s^-[h_s(1)],  \\
& u_+=\eta_s^\vee v_+, \ u_-=\eta_s^\vee u_- \in \mathcal{Y}_s^-[h_s(0)]. 
\end{align*}
Then we have 
$$W^+(0)v_+=0$$
Let be $[E] \in \mbox{Ext}^1(\mathcal{Y}_s^-, \mathcal{X}_s^-)$, then we have an exact sequence of $W(p)$-module,
$$0 \longrightarrow E_0 \longrightarrow E \overset{\pi}{\longrightarrow} \mathcal{Y}_s^- \longrightarrow 0,$$
where $E_0$ is isomorphic to $\mathcal{X}_s^-$.  Fix elements of $E_0$ by the following way,
$$v_+^{(1)}, \ v_-^{(1)}=W^-(0)v_+^{(1)} \in E_0[h_s(1)] \simeq \mathbb C^2.$$
Then we have $W^+(0)v_+^{(1)}=0$. And take elements of $E$ by the following way,
\begin{align*}
& \tilde{v}_+ \in E[h_s(1)] \rightarrow v_+ \in \mathcal{Y}_s^-, \\
& \tilde{v}_-=W^-(0)\tilde{v}_+, \\
& \tilde{u}_\pm=\eta_s^\vee \tilde{v}_\pm
\end{align*}
Then we have $\pi(\tilde{v}_\pm)=u_\pm$ and $W^+(0)\tilde{v}_+=0$.
\par
For $W(p)$-module $M \in W(p)$-mod, we define $Q : M \rightarrow M$ by the following way,
$$Q |_{M[h]}=(T(0)-h) \ \mbox{id}.$$
Then the $\mathbb C$ linear map $Q : M \rightarrow M$ is $W(p)$-module map, and satisfies $Q^n=0$ for some $n \geqq 1$.  Then we have a commutative diagram,
\begin{equation}\tag{5.15}
\begin{CD}
0 @>>> E_0 @>>> E @>>> \mathcal{Y}_s^- @>>> 0 \\
@.   @VV{Q}V   @VV{Q}V   @VV{Q}V   @. \\
0 @>>> E_0 @>>> E @>>> \mathcal{Y}_s^- @>>> 0.
\end{CD}
\end{equation}
Since the map $Q=0$ on $\mathcal{Y}_s^-$ and on $E_0$, we have
$$Q(E) \subseteq E_0 \subseteq E$$
and $Q^2=0$ on $E$.  Therefore $Q$ factor through
$$Q : E \overset{\bar{\pi}}{\longrightarrow} \mathcal{X}_s^- \overset{\bar{Q}}{\longrightarrow} E_0$$
where $\bar{\pi} : E \overset{\pi}{\rightarrow} \mathcal{Y}_s^- \rightarrow \mathcal{X}_s^-$.  Since $Q$ is $W(p)$-module map their exact constant $\gamma$, such that
\begin{equation}\tag{5.16}
Q(\tilde{v}_\pm) = \gamma \ \tilde{v}_\pm^{(1)}.
\end{equation}
We show that if $\gamma=0$, then $[E]=0$ in $\mbox{Ext}^1(\mathcal{Y}_s^-, \mathcal{X}_s^-)$.  Consider the exact sequence
$$0 \longrightarrow (\mathcal{X}_s^+)^2 \longrightarrow \mathcal{Y}_s^- \longrightarrow \mathcal{X}_s^- \longrightarrow 0,$$
then we have an exact sequence
$$0 \longrightarrow \mbox{Ext}^1(\mathcal{Y}_s^-, \mathcal{X}_s^-) \longrightarrow \mbox{Ext}^1((\mathcal{X}_s^+)^2, \mathcal{X}_s^-).$$
Therefore to prove $[E]=0$, it is sufficient to prove that
$$\eta_s(\tilde{u}_\pm)=0.$$
By Proposition 2-12, (2.40) we see that
\begin{align*}
\eta_s(\tilde{u}_\pm) & = \eta_s\eta_s^\vee (\tilde{v}_\pm) \\
& = c(T(0)-h_s(1))\tilde{v}_\pm,
\end{align*}
for some $c \ne 0$.  By the assumption $\gamma=0$, we have $\eta_s(\tilde{u}_\pm)=0$.
\par
This show that
$$\dim_{\mathbb C} \mbox{Ext}^1(\mathcal{Y}_s^-, \mathcal{X}_s^-) \leqq 1.$$
\end{proof}

\vspace{0.2in}
\noindent
{\bf Proposition 5-17.} 
$$\mbox{Ext}^1(\mathcal{P}_s^-, \mathcal{X}_s^-) = 0.$$

\begin{proof}
Consider the exact sequence
$$0 \longrightarrow \mathcal{X}_s^- \longrightarrow \mathcal{P}_s^- \overset{\pi}{\longrightarrow} \mathcal{Y}_s^- \longrightarrow 0.$$
Then we have an exact sequence
$$0 \longrightarrow \mbox{Hom}(\mathcal{X}_s^-, \mathcal{X}_s^-) \longrightarrow \mbox{Ext}^1(\mathcal{Y}_s^-, \mathcal{X}_s^-) \longrightarrow \mbox{Ext}^1(\mathcal{P}_s^-, \mathcal{X}_s^-) \longrightarrow 0.$$
By Proposition 5-16, we have $\mbox{Ext}^1(\mathcal{P}_s^-, \mathcal{X}_s^-)=0$.
\end{proof}

\vspace{0.2in}
\noindent
{\bf Proposition 5-18.} 
$\mathcal{P}_s^-$ is projective cover of simple $W(p)$-module $\mathcal{X}_s^-$.

\begin{proof}
By Proposition 5-15 and Proposition 5-17, the $W(p)$-module $\mathcal{P}_s^-$ is projective.
\end{proof}

\vspace{0.2in}
The following propositions can be easily proved by using the above propositions.

\vspace{0.2in}
\noindent
{\bf Proposition 5-19.} 
We have
$$\mbox{Ext}^1(\mathcal{V}_s^-, \mathcal{X}_s^-)=0,$$
$$\mbox{Ext}^1(\mathcal{V}_s^{-\vee}, \mathcal{X}_s^-)=0.$$

\vspace{0.2in}
\noindent
{\bf Proposition 5-20.}
\begin{enumerate}
\item $\mbox{Ext}^1(\mathcal{V}_s^-, \mathcal{V}_s^+) \simeq \mathbb C$, $\mbox{Ext}^1(\mathcal{V}_s^{-\vee}, \mathcal{V}_s^{+\vee}) \simeq \mathbb C$ and these two vector spaces are generated by $[\mathcal{P}_s^-] .$
\item $\mbox{Ext}^1(\mathcal{V}_s^-, \mathcal{V}_s^{+\vee}) =0$, $\mbox{Ext}^1(\mathcal{V}_s^{-\vee}, \mathcal{V}_s^+) =0$.
\end{enumerate}

\vspace{0.2in}
\noindent
{\bf Proposition 5-21.} 
We have
\begin{enumerate}
\item $D(\mathcal{P}_s^-)\simeq \mathcal{P}_s^-$.
\item $\mathcal{P}_s^-$ is an injective module.
\end{enumerate}

\vspace{0.5in}
\section{Category equivqlent of $W(p)$-mod and $\bar U_q(sl_2)$-mod}

In this section we prove that two abelian categories $W(p)$-mod and $\bar U_q(sl_2)$-mod are equivalent as abelian categories.  This is conjected in [FGST1], [FGST2]

\vspace{0.2in}
\subsection{Quantum group $\bar U_q(sl_2)$}
We fix positive integer $p \geqq 2$, and set $q=e^{\pi i / p}$.  We introduce the restricted quantum group $\bar U_q(sl_2)=\bar U$.  We will follow the articles of Feigin et al.  [FGST1], [FGST2] and Kondo and Saito [KoS].
\par
For each integer $n$, we set
\begin{equation}\tag{6.1}
[n]=\frac{q^n-q^{-n}}{q-q^{-1}}
\end{equation}
The restricted quantum group $\bar U_q(sl_2)$ is an associative $\mathbb C$-algebra with the unit, which is generated by $E, F, K, K^{-1}$ satisfying the following fundamental relations
\begin{equation}\tag{6.2}
KK^{-1}=K^{-1}K=1, \ KEK^{-1}=q^2E, \ KFK^{-1}=q^{-2}F
\end{equation}
$$EF-FE=\frac{K-K^{-1}}{q-q^{-1}}$$
$$K^{2p}=1, \ E^p=F^p=0$$
This is finite dimensional $\mathbb C$-algebra, and has a structure of Hopf-algebra.
\par
Let $\bar U$-mod be the abelian category of finite dimensional $\bar U$-modules.  Then it is known in [FGST1], [FGST2], [KoS].

\vspace{0.2in}
\noindent
{\bf Proposition 6-1.}  
The abelian category $\bar U$-mod has blocks decomposition
$$\bar U\mbox{-mod}=\sum_{s=0}^p C_s(\bar U)$$
where $C_0(\bar U)$ and $C_p(\bar U)$ are semi-simple categories whose have only one simple object, respectively.   For $1 \leqq s \leqq p-1$, $C_s(\bar U)$ are all isomorphic each other as abelian categories.  The category $C_s(\bar U)$ is Artinian and Neotherian, and the number of simple object is two.  We denote this abelian category $C(\bar U)$, and denote simple object $\mathcal{X}^\pm(\bar U)=\mathcal{X}^\pm$ and their projective cover $\mathcal{P}^\pm(\bar U)=\mathcal{P}^\pm$.  Set $P(\bar U)=P^+(\bar U)\oplus P^-(\bar U)=P^+ \oplus P^-$.
\par
And consider finite dimensional $\mathbb C$-algebra
\begin{equation}\tag{6.3}
B(\bar U)=\mbox{End}_{C(\bar U)}(\bar P(\bar U)).
\end{equation}
The following proposition is known in [FGST1], [FGST2], [KoS].

\vspace{0.2in}
\noindent
{\bf Theorem 6-2.} 
$B(\bar U)$ is 8 dimensional algebra of the form;
\begin{align}\tag{6.4}
B(\bar U)= 
& \mbox{End}_C(\mathcal{P}_s^+, \mathcal{P}_s^+) \oplus \mbox{End}_C(\mathcal{P}_s^-, \mathcal{P}_s^-) \\
& \qquad \oplus \mbox{Hom}_C(\mathcal{P}_s^+, \mathcal{P}_s^-) \oplus \mbox{Hom}_C(\mathcal{P}_s^-, \mathcal{P}_s^+), \notag
\end{align}
and generated by
\begin{equation}\tag{6.5}
\mbox{Hom}_C(\mathcal{P}_s^\pm, \mathcal{P}_s^\mp) = C\tau_1^\pm \oplus C\tau_2^\pm,
\end{equation}
with the relations;
\begin{align}\tag{6.6}
& \tau_i^\pm \tau_i^\mp=0, \quad i=1, 2, \\
& \tau_1^\pm \tau_2^\mp=\tau_2^\pm \tau_1^\mp, \notag
\end{align}

\vspace{0.2in}
Now we consider a $\mathbb C$-linear abelian category $C$ with the following properties;
\begin{enumerate}
\item $C$ is Neotherian and Artinian.
\item The set of equivalence classes of simple objects is finite, say $\{ S_1, \dots , S_N \}$.
\end{enumerate}
\par
Denote the projective cover of $S_i$ by $P_i$.  And set $P=\sum_{i=1}^N P_i$.  Consider the Endmorphism algebra of $\mathcal{P}$,
$$B(C)=\mbox{End}_C(P)$$
Then $B(C)$ is a finite dimensional algebra over $\mathbb C$.
\par
Denote by mod $B(C)$, the abelian category of finite dimensional right $B(C)$-modules.  Then the following proposition is well known.

\vspace{0.2in}
\noindent
{\bf Proposition 6-3.} 
\begin{align*}
\Phi : C & \longrightarrow \mbox{mod-} B(C) \\
M & \longrightarrow \mbox{Hom}_{C}(P, M)
\end{align*}
is equivalence of abelian categories.

\vspace{0.5in}
\subsection{Categorical equivalence of two abelian category $W(p)$-mod and $\bar U$-mod}
We showed that the abelian category $W(p)$-mod has the block decomposition
$$W(p)\mbox{-mod}=\sum_{s=0}^p C_s,$$
and that $C_0$ and $C_p$ are semi-simple categories whose simple objects are $\mathcal{X}_0$ and $\mathcal{X}_p$, respectively.  On the other hand for $1\leqq s \leqq p-1$, each abelian category $C_s$ has two simple objects $\mathcal{X}_s^+$ and $\mathcal{X}_s^-$.

\vspace{0.2in}
Now for each $1\leqq s \leqq p-1$, consider $\mathcal{P}_s=\mathcal{P}_s^+ \oplus \mathcal{P}_s^- \in C_s$, and define finite dimensional $\mathbb C$-algebra $B_s$ as follows
$$B_s=\mbox{End}_{\mathbb C}(C_s).$$

\vspace{0.2in}
\noindent
{\bf Theorem 6-4.} 
For each $1\leqq s \leqq p-1$, the finite dimensional algebra $B_s$ is isomorphic to $B(\bar U)$

\begin{proof}
By proposition 4-5, it is easy to show that $B_s$ is isomorphic to $B(\bar U)$ as an algebra over $\mathbb C$.
\end{proof}

\vspace{0.2in}
\par
So by Proposition 6-3, we have the following main theorem of this section.

\vspace{0.2in}
\noindent
{\bf Theorem 6-5.}
Two abelian categories $W(p)$-mod and $\bar U_q(sl_2)$-mod are equivalent.

\vspace{0.5in}
\subsection{Length of the Jordan blocks}
For each $M\in W(p)$-mod, we define $l(M)\in \mathbb Z_{\geqq 0}$ by $l(M)=$ max $\{ n \in \mathbb Z_{\geqq 0}$ ; $(T(0)-h)^n v \ne 0$ for some $h \in \mathbb C$, $v \in M[h] \}$.
\par
Then we obtain the following proposition.

\vspace{0.2in}
\noindent
{\bf Proposition 6-6.} 
\begin{enumerate}
\item For each $M \in W(p)$-mod we have $l(M) \leqq 1$
\item Any indecomposable module $M$ in $W(p)$-mod such that $l(M)=1$ is equivalent to $M \simeq \mathcal{P}_s^+$ for some $s$ such that $1\leqq s \leqq p-1$.
\end{enumerate}

\end{document}